\documentclass{article}%
\usepackage{epsfig}
\usepackage{enumerate,graphicx}
\usepackage{amsfonts}
\usepackage{graphicx}
\usepackage{amsmath,amsthm,amssymb}
\usepackage{amssymb}
\usepackage{appendix}
\DeclareGraphicsExtensions{.pdf, .jpg, .png , .bmp, .ps}
\usepackage{tikz}
\providecommand{\U}[1]{\protect\rule{.1in}{.1in}}
\newtheorem{thm}{Theorem}[section]
\newtheorem{lem}{Lemma}[section]
\newtheorem{prop}{Proposition}[section]
\newtheorem{cor}{Corollary}[section]
\newtheorem{Defi}{Definition}[section]

\theoremstyle{remark}
\newtheorem{rem}[thm]{\bf Remark}
\newtheorem{defn}[thm]{Definition}

\newcommand{\tq}{\,:\,}\newcommand\E{{\mathbb {E}}}
\newcommand\F{{\mathcal {F}}}
\newcommand\I{{\bf 1}}

\newcommand\bkE{{\mathbb {E}}}
\newcommand\p{{\mathbb {P}}}

\DeclareMathOperator{\Var}{V}

\def\1{{{\mbox{${\rm{1\negthinspace\negthinspace I}}$}}}}

\newcommand{\eref}[1]{(\ref{#1})}

\newcommand\beq{\begin{equation}}
\newcommand\eeq{\end{equation}}

\begin{document}

\title{Behavior of the Wasserstein distance between
the empirical and the marginal distributions  of  stationary $\alpha$-dependent sequences}

\author{J\'er\^ome Dedecker\footnote{Universit\'e Paris Descartes,
Sorbonne Paris Cit\'e,
Laboratoire MAP5
and CNRS UMR 8145.},
Florence Merlev\`ede \footnote{Universit\'{e} Paris Est,  UPEM, UPEC, LAMA and CNRS UMR 8050.}
}

\date{}

\maketitle

\abstract{We study the Wasserstein distance of order 1 between the empirical distribution and the 
marginal distribution of  stationary  $\alpha$-dependent sequences. We prove some moments inequalities 
of order $p$ for any $p\geq 1$, and we give some conditions under which the central limit theorem 
holds. We apply our results to unbounded functions of expanding maps of the interval with a neutral fixed point at zero. The moment inequalities for the Wasserstein distance are similar to the well known 
von Bahr-Esseen  or Rosenthal bounds for partial sums, and seem to be new even in the 
case of independent and identically distributed random variables.

\medskip

\noindent{\bf Running head.} Empirical Wasserstein 
distances for dependent sequences. 

\medskip


\medskip


\section{Introduction}
Let $(X_i)_{i \in {\mathbb Z}}$ be a stationary sequence of integrable real-valued random 
variables, with common marginal distribution $\mu$. Let 
$\mu_n$ be the empirical measure of $\{X_1, \ldots , X_n\}$, that is
$$
  \mu_n=\frac 1 n \sum_{k=1}^n \delta_{X_k} \, .
$$
In this paper, we study the behavior of the quantity $W_1(\mu_n, \mu)$ for a large class
of stationary sequences, where $W_1(\mu_1, \mu_2)$ is the {\it Wasserstein distance} of order 1
between two probability measures $\mu_1, \mu_2$ having finite first moments. The precise definition is as follows:
\begin{equation}\label{start}
   W_1(\mu_1, \mu_2)= \inf_{\pi \in M(\mu_1, \mu_2)} \int |x-y| \pi(dx,dy) \, ,
\end{equation}
where $M(\mu_1, \mu_2)$ is the set of probability measures on ${\mathbb R}^2$ with marginal
distributions $\mu_1$ and $\mu_2$. The distance $W_1$ belongs to the general class of minimal distances,
as the total variation distance. Since  the cost function $c_1(x,y)=|x-y|$ is regular, $W_1$ can be used to compare two singular measures, which is not possible with the total variation distance, whose 
cost function is given by the discrete metric $c_0(x,y)={\bf 1}_{x \neq y}$. 

  The quantity $W_1(\mu_n, \mu)$ appears very frequently in statistics, and can be understood from many
  points of view:
 \begin{itemize} 
  \item The well known dual representation of $W_1$ implies that
\begin{equation}\label{dual}
  W_1(\mu_n, \mu)= \sup_{f \in \Lambda_1} \left | \frac 1 n \sum_{k=1}^n \left ( f(X_k) - \mu(f)
  \right ) \right| \, ,
\end{equation}
where $\Lambda_1$ is the set of  Lipschitz functions $f$ from 
${\mathbb R}$ to ${\mathbb R}$
such that $|f(x)-f(y)|\leq |x-y|$. Hence, 
$W_1(\mu_n, \mu)$ is a measure of the concentration of $\mu_n$ around $\mu$ through the class $\Lambda_1$. 
  \item In the one dimensional setting the minimization problem \eref{start} can be explicitely solved,
  and leads to the expression
\begin{equation}\label{empquant}
  W_1(\mu_n, \mu)= \int_0^1 |F_n^{-1}(t)-F^{-1}(t)| dt \, ,
\end{equation}
where $F_n$ and $F$ are the distribution functions of $\mu_n$ and $\mu$, and $F_n^{-1}$ 
and $F^{-1}$ are their usual generalized inverses. Hence $W_1(\mu_n, \mu)$ is the ${\mathbb L}^1$-distance 
between the empirical quantile function $F_n^{-1}$ and the quantile function  of $\mu$. 
\item Starting from \eref{empquant}, it follows immediately that
\begin{equation}\label{empL1}
W_1(\mu_n, \mu)= \int_{\mathbb R} |F_n(t)-F(t)| dt \, .
\end{equation}
Hence $W_1(\mu_n, \mu)$ is the ${\mathbb L}^1$-distance 
between the empirical distribution function $F_n$ and the distribution function of $\mu$.
\end{itemize}
At this point, it should be clearly quoted that, if \eref{empquant} and \eref{empL1}
have no analogue in higher dimension, the dual expression 
\eref{dual} is very general
and holds if the $X_i$'s take their values in a Polish space
${\mathcal X}$, as soon as the cost function $c$
is a lower semi-continuous metric (the class $\Lambda_1$ being the class of 1-Lipschitz functions from ${\mathcal X}$ to ${\mathbb 
 R}$ with
respect to $c$). 
  
  Assume now that the sequence $(X_i)_{i \in {\mathbb Z}}$ is ergodic. Since $\mu$ has a finite first moment, it is well known that
$W_1(\mu_n, \mu)$ converges to zero almost surely, and that ${\mathbb E}(W_1(\mu_n, \mu))$ converges to zero (this is a uniform version of Birkhoff's ergodic theorem, 
which can be easily deduced from the Glivenko-Cantelli theorem for ergodic sequences). However, without additional asumptions
on $\mu$ the rate of convergence can be arbitrarily slow. 
 
   The purpose of this paper is to give some conditions under which the central limit theorem (CLT) holds
   (meaning that $\sqrt n W_1(\mu_n, \mu)$ converges in distribution to a certain law), and to prove some 
   inequalities for $\|W_1(\mu_n, \mu)\|_p$ when $p \geq 1$ (von Bahr-Esseen type inequalities for 
   $p \in (1,2)$ and Rosenthal type inequalities for $p>2$). We will do this for the class of 
   $\alpha$-dependent sequences, which is quite natural in this context, since the related dependency coefficients
   are defined through indicator of half lines. Hence our results apply to mixing sequences in 
   the sense of Rosenblatt \cite{R}, but also to many other dependent sequences including a large class
   of one dimensional dynamical systems. We shall illustrate our results through the examples of
   Generalized Pomeau-Manneville maps, as defined in \cite{DGM}. 
  
The central limit question for $\sqrt n W_1(\mu_n, \mu)$ has been already investigated 
for dependent sequences in the papers by D\'ed\'e \cite{D} and Cuny \cite{C} (see Sections \ref{CLTL1} and
\ref{quantile} for more details). This is not the case of the
upper bounds for $\|W_1(\mu_n, \mu)\|_p$, even for  sequences of independent and identically distributed
(i.i.d.) random variables (except for $p=1$, see for instance \cite{BL}). Hence, for $p>1$,
our moment bounds seem to be new  even in the i.i.d. context.

Thanks to the relation \eref{empL1}, the central limit question for $\sqrt n W_1(\mu_n, \mu)$ is closely 
related to the empirical central limit theorem in ${\mathbb L}^1(dt)$, as first quoted by 
del Barrio, Gin\'e and Matr\'an \cite{BGM}. We shall deal with the more general  central limit question 
for ${\mathbb L}^1(m)$-valued random variables in the separate Section \ref{CLTL1}. In Section
\ref{quantile}, 
we shall express some of our conditions in terms of the quantile function of $X_0$, in the spirit of 
Doukhan, Massart and Rio \cite{DMR}. It will then be easier  to compare our conditions for the CLT to 
previous ones
in the literature.

For $r>1$, the quantity $W_r^r(\mu_n, \mu)$ may be defined as in \eref{start},  with the cost function 
$c_r(x,y)=|x-y|^r$ instead of $c_1$ ($W_r$ is the Wasserstein distance of order $r$).
In the i.i.d. case,  some sharp upper bounds on 
${\mathbb E}(W^r_r(\mu_n, \mu))$ are given in the recent
paper \cite{BL}. In particular, if $\mu$ has an absolutely component  with respect
to the Lebesgue measure which does not vanishes on the support of $\mu$, then the optimal rate
$n^{-r/2}$ can be reached. But in general, the rate can be much slower. Note that for $W^r_r(\mu_n, \mu)$
there is no such nice dual expression as  \eref{dual}. However the minimization problem can still be 
explicitely solved and implies that $W_r(\mu_n, \mu)$ is the ${\mathbb L}^r$-distance between 
$F_n^{-1}$ and $F^{-1}$. There is no simple way to express $W^r_r(\mu_n, \mu)$ in terms of 
$F_n$ and $F$ (as in \eref{empL1}), but the following upper bound due to 
\`Ebralidze \cite{E} holds: 
\begin{equation}\label{Ebr}
W_r^r(\mu_n, \mu) \leq \kappa_r \int_{\mathbb R}
 |x|^{r-1} |F_n (x) -F (x) | dx \, ,
\end{equation}
where $\kappa_r= 2^{r-1}r$. Starting from this inequality,  we shall also give some upper 
bounds on $\|W^r_r(\mu_n, \mu)\|_p$ for $p\geq 1$, but it is very likely that these bounds
can be improved by assuming the existence of an absolutely regular component for $\mu$, as 
in \cite{BL}. 

To be complete, let us mention the recent paper by Fournier and Guillin \cite{FG}, who give some upper bounds
for ${\mathbb E}(W^r_r(\mu_n, \mu))$ in any dimension, starting from an inequality which can be viewed as 
a $d$-dimensional analogue of \eref{Ebr}. Note that the case of $\rho$-mixing sequences is also considered 
in this paper. 

\section{Definitions and notations} \label{DefNot}
\setcounter{equation}{0}

In this section, we give the notations and definitions which we will used all
along the paper.  

Let us start with the  notation
 $a_n \ll b_n$, which means that there exists a numerical constant $C$ not
depending on $n$ such that  $a_n \leq  Cb_n$, for all positive integers $n$.
\subsection{Stationary sequences and dependency coefficients}

Let  $(\Omega ,\mathcal{A}, \p)$ be a probability space, and 
$T :\Omega \mapsto \Omega $ be a bijective bi-measurable
transformation preserving the probability ${\p}$. Let
${\mathcal F}_0$ be a sub-$\sigma$-algebra of $\mathcal{A}$
satisfying ${\mathcal F}_0 \subseteq T^{-1}({\mathcal
F}_0)$. We say that the couple $(T, \p)$ is ergodic if any 
$A \in {\mathcal A}$ satisfying $T(A)=A$ has probability 0 or 1. 

Let $X_0$ be an ${\mathcal F}_0$-measurable and integrable real-valued random variable
with distribution $\mu$. 
Define the stationary sequence 
 ${\bf X}=(X_i)_{i \in {\mathbb Z}}$ by $X_i=X_0 \circ T^i$.
 
Let us first define the tail and quantile functions of the 
random variable $X_0$. 

\medskip

\begin{defn}
The tail function  $H:{\mathbb R}^+ 
\rightarrow  [0,1]$ of $X_0$ is defined by 
$H(t)= {\mathbb P}(|X_0|>t)$.  The quantile function 
$Q: [0,1]  \rightarrow  {\mathbb R}^+$ of $X_0$
is the generalized  inverse 
of $H$, that is 
$$
Q(u)= \inf \left \{ t\geq 0: f(t) \leq u \right \} \, .
$$
\end{defn}

\medskip


Let us now  define the dependency coefficients of the sequence 
$(X_i)_{i \in {\mathbb Z}}$. These coefficients
are less restrictive than the usual mixing coefficients of Rosenblatt \cite{R}. 

\medskip

\begin{defn}
For any integrable random variable $Z$, let 
$Z^{(0)}=Y- \E(Z)$.
For any random variable $Y=(Y_1, \cdots, Y_k)$ with values in
${\mathbb R}^k$ and any $\sigma$-algebra ${\mathcal F}$, let
\[
\alpha({\mathcal F}, Y)= \sup_{(x_1, \ldots , x_k) \in {\mathbb R}^k}
\left \| \E \left(\prod_{j=1}^k (\I_{Y_j \leq x_j})^{(0)} \Big | {\mathcal F} \right)
-\E \left(\prod_{j=1}^k (\I_{Y_j \leq x_j})^{(0)}  \right) \right\|_1.
\]
For the stationary sequence ${\bf X}=(X_i)_{i \in {\mathbb Z}}$, let 
\begin{equation}
\label{defalpha} \alpha_{k, {\bf X}}(n) = \max_{1 \leq l \leq
k} \ \sup_{ n\leq i_1\leq \ldots \leq i_l} \alpha({\mathcal F}_0,
(X_{i_1}, \ldots, X_{i_l})).
\end{equation}
Note that $\alpha_{1, {\bf X}}(n)$ is then simply given by 
\begin{equation}
\alpha_{1, {\bf X}}(n)= \sup_{x \in {\mathbb R}}\left \|\E\left (\I_{X_n \leq x}|{\mathcal F}_0\right ) - F(x)
\right \|_1 \, ,
\end{equation} 
where $F$ is the distribution function of $\mu$. 
\end{defn}

\medskip

All the results of Section \ref{Sec:main} below involve only the coefficients 
$\alpha_{1, {\bf X}}(n)$, except for the Rosenthal bounds (Subsection \ref{Sec:Ros}) for which the
coefficient $\alpha_{2, {\bf X}}(n)$ is needed.

\subsection{Intermittent maps} 
Let us first recall the definition of the generalized Pomeau-Manneville maps introduced in 
\cite{DGM}. 

\begin{defn}
A map $\theta:[0,1] \to [0,1]$ is a generalized Pomeau-Manneville
map (or GPM map) of parameter $\gamma \in (0,1)$ if there exist
$0=y_0<y_1<\dots<y_d=1$ such that, writing $I_k=(y_k,y_{k+1})$,
\begin{enumerate}
\item The restriction of $\theta$ to $I_k$ admits a $C^1$ extension
$\theta_{(k)}$ to $\overline{I_k}$.
\item For $k\geq 1$, $\theta_{(k)}$ is $C^2$ on $\overline{I_k}$, and $|\theta_{(k)}'|>1$.
\item $\theta_{(0)}$ is $C^2$ on $(0, y_1]$, with $\theta_{(0)}'(x)>1$ for $x\in
(0,y_1]$, $\theta_{(0)}'(0)=1$ and $\theta_{(0)}''(x) \sim c
x^{\gamma-1}$ when $x\to 0$, for some $c>0$.
\item $\theta$ is topologically transitive.
\end{enumerate}
\end{defn}
The third condition ensures that $0$ is a neutral fixed point
of $\theta$, with $\theta(x)=x+c' x^{1+\gamma} (1+o(1))$ when $x\to 0$.
The fourth condition is necessary to avoid situations where
there are several absolutely continuous invariant measures, or
where the neutral fixed point does not belong to the support of
the absolutely continuous invariant measure.

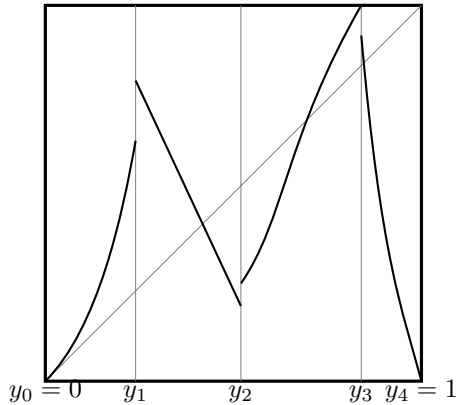
\begin{figure}[htb]
\centering
  \begin{tikzpicture}[scale=1]
  \draw[very thick] (0,0) rectangle (5,5);
  \draw[gray, very thin]
      (1.2,0) -- +(0, 5)
      (2.6,0) -- +(0, 5)
      (4.2,0) -- +(0, 5)
      (0,0)   -- (5,5);
  \draw[thick]
      (0,0) .. controls +(45:1) and +(-100:1) .. (1.2, 3.2)
      (1.2, 4) -- (2.6, 1)
      (2.6, 1.3) .. controls +(55:1) and +(-120:2) .. (4.2, 5)
      (4.2,4.6) .. controls +(-85:3) and +(105:1) .. (5, 0);
  \foreach \x/\ytext in {0/$y_0=0$, 1.2/$y_1$, 2.6/$y_2$, 4.2/$y_3$, 5/$y_4=1$}
      \node[above] at (\x, -0.4) {\ytext};
  \end{tikzpicture}
\caption{The graph of a GPM map, with $d=4$}
\end{figure}
The following  well known example of GPM map with only two branches 
has been introduced 
by  Liverani, Saussol and Vaienti \cite{LSV}:
\beq \label{LSVmap}
   \theta(x)=
  \begin{cases}
  x(1+ 2^\gamma x^\gamma) \quad  \text{ if $x \in [0, 1/2[$}\\
  2x-1 \quad \quad \quad \ \  \text{if $x \in [1/2, 1]$.}
  \end{cases}
\eeq

As quoted in \cite{DGM}, a GPM map $\theta$ admits 
a unique invariant absolutely continuous (with respect to the Lebesgue measure) probability  $\nu$ with density $h$. 
Moreover, it is ergodic, has full support, and $x^\gamma h(x)$ is bounded from above and below.

We shall illustrate each result of Section \ref{Sec:main} by controlling, on the probability space $([0,1], \nu)$, the
quantity 
$
W_1(\tilde \mu_n, \mu)
$, 
where 
\beq \label{empGPM}
 \tilde \mu_n= \frac 1 n \sum_{k=1}^n 
 \delta_{g \circ \theta^k}  \, ,
\eeq
$\theta$ is a GPM map, $g$ is a monotonic function from $(0,1)$ to ${\mathbb R}$ (which can blow up near
$0$ or $1$), and $\mu$ is the distribution of $g$. 

To do this, we  go back to the Markov chain associated to $\theta$, as we describe now. 
Let first $K$ be the
 Perron-Frobenius operator of $\theta$ with respect to $\nu$, defined as follows:
for any  functions $u, v$ in ${\mathbb L}^2([0,1], \nu)$
\begin{equation}\label{Perron}
\nu(u \cdot  v\circ \theta)=\nu(K(u) \cdot  v) \, .
\end{equation}
The relation \eref{Perron} states  that $K$ is  the adjoint operator of the isometry $U: u \mapsto u\circ \theta$
acting on ${\mathbb L}^2([0,1], \nu)$. It is easy to see that the operator $K$ is a transition kernel, 
and that $\nu$ is invariant by $K$. 
Let now 
${\bf Y}=(Y_i)_{i\geq 0}$
be a stationary Markov chain with invariant measure $\nu$
and transition kernel $K$. It is well known (see for instance Lemma XI.3 in
\cite{HH}) that on the probability space $([0, 1], \nu)$, the
random vector $(\theta, \theta^2, \ldots , \theta^n)$ is distributed as
$(Y_n,Y_{n-1}, \ldots, Y_1)$.

Let $T$ be the shift operator from ${[0,1]}^{\mathbb
Z}$ to ${[0,1]}^{\mathbb Z}$ defined by
$(T(x))_i=x_{i+1}$, and let $\pi_i$ be the projection from
${[0,1]}^{\mathbb Z}$ to ${[0,1]}$ defined by
$\pi_i(x)=x_i$.  By Kolmogorov's extension theorem,
there exists a shift-invariant probability ${\mathbb P}$ on
$({[0,1]}^{\mathbb Z}, ({\mathcal B}([0,1]))^{\mathbb
Z})$, such that $\pi=(\pi_i)_{i\geq 0}$ is distributed as ${\bf
Y}$. 

Let then $X_0=g \circ \pi_0$ and $X_i=X_0 \circ T^i= g \circ \pi_i$, 
and define 
${\mathcal F}_0= \sigma (\pi_i, i \leq 0)$. From the above considerations, we infer 
that the two random variables 
 $W_1(\mu_n, \mu)$ (defined on the probability space $({\mathbb R}^{\mathbb Z}, {\mathbb P})$) and $W_1(\tilde \mu_n, \mu)$ (defined on the probability space $([0,1], \nu)$) 
have the same distribution. Hence, any information on the distribution of
$W_1(\tilde \mu_n, \mu)$ can be derived from the distribution of 
$W_1(\mu_n, \mu)$.


From Proposition 1.17 (and the comments right after)  in \cite{DGM}, we know that
for any positive integer $k$, there exist two  positive constants $C$ and $D$ such
that, for any $n>0$,
\[
 \frac{D}{n^{(1-\gamma)/\gamma}} \leq \alpha_{k,{\pi}}(n) \leq
  \frac{C}{n^{(1-\gamma)/\gamma}} \, .
\]
Since $X_i= g\circ \pi_i$, and since $g$ is monotonic, it follows  immediately that
\begin{equation}\label{boundalpha}
\alpha_{k,{\bf X}}(n) \leq \alpha_{k,{\pi}}(n) \leq
  \frac{C}{n^{(1-\gamma)/\gamma}} \, .
\end{equation}
This control of the coefficients $\alpha_{k,{\bf X}}(n)$
(for $k=1$ or $k=2$) and a control 
of the tail $\nu(|g|>t)$ are all  we need to apply the results of Section \ref{Sec:main}
to the random variable $W_1(\tilde \mu_n, \mu)$.

\section{CLT and moment bounds}\label{Sec:main}
\setcounter{equation}{0}
In all this section, we use the notations of Sections 1 and 2.
\subsection{Central limit theorem}\label{CLT}
Our first result is a central limit theorem for $W_1(\mu_n, \mu)$. It is a straightforward consequence of a 
 CLT in ${\mathbb L}^1(m)$  for the empirical distribution function given in Proposition \ref{prop:alpha} of Subsection \ref{intro} (it  suffices 
 to consider the case where $m$ is the Lebesgue measure
 on ${\mathbb R}$ and to use the continuous mapping theorem).
\begin{prop}\label{prop:alphaW1}
Assume that the couple $(T, {\mathbb P})$ is ergodic,
 and  that 
\begin{equation}\label{condalpha}
\int_0^{\infty} \sqrt{ \sum_{k=0}^\infty \min \left\{\alpha_{1, \bf X}(k), H(t)
\right \} } \ dt < \infty \, .
\end{equation}
Then $\sqrt n W_1(\mu_n, \mu)$ converges in distribution to the random variable
$\int |G(t)| \ dt$, where $G$ is a Gaussian random variable in ${\mathbb L}^1(dt)$ whose covariance 
function may be described as follows:
for any $f,g$ in ${\mathbb L}_\infty(\mu)$,
\begin{equation}\label{defcovW1}
\mathrm{Cov}\left(\int f(t) G(t) dt, \int g(t)G(t) dt \right)= \sum_{k \in {\mathbb Z}} {\mathbb{E}}
\left( \iint f(t)g(s) ({\mathbf{1}}_{X_0 \leq t}-F(t))   ({\mathbf{1}}_{X_k \leq s}-F(s)) \ dt  ds
 \right)\, .
\end{equation}
\end{prop}

\begin{rem} Let $m$ be a nonnegative integer. 
As usual, the stationary sequence $\bf X$ is 
$m$-dependent if 
$\sigma(X_i, i \leq 0)$ is independent of 
$\sigma(X_i, i \geq m+1)$,  and $m=0$ 
corresponds to  the i.i.d. case. 
In the $m$-dependent case,  the condition \eref{condalpha} 
becomes simply
\begin{equation}\label{condiid}
\int_0^{\infty} \sqrt{  H(t) } \ dt < \infty \, , 
\end{equation}
which is exactly the condition given by del Barrio, Gin\'e and Matr\'an \cite{BGM} in the i.i.d. case. Note that these authors also 
proved that, in the i.i.d. case, the condition 
\eref{condiid} is necessary and sufficient for the stochastic boundedness of $\sqrt n W_1(\mu_n, \mu)$.

In the dependent context, other general criteria have 
been proposed by D\'ed\'e \cite{D} and Cuny \cite{C}. We shall
discuss these conditions in Sections \ref{CLTL1} and
 \ref{quantile}, and show that, in 
the $\alpha$-dependent case, the condition \eref{condalpha} is weaker than the corresponding 
condition
obtained by applying the criteria by D\'ed\'e or Cuny.
\end{rem}

\noindent{\bf Example.}
Let $\theta$ be a GPM map of parameter $\gamma \in (0,1/2)$, with absolutely 
continuous invariant probability $\nu$. Let $\tilde \mu_n$ be defined as in \eref{empGPM}, where
$g$ is a monotonic function from $(0,1)$ to ${\mathbb R}$. 
Let then $(X_i)_{i \in {\mathbb Z}}$ be the stationary sequence constructed in Subsection 2.2, whose 
dependency coefficients $\alpha_{k, {\bf X}}(n)$ satisfy \eref{boundalpha}. Note that
$H(t)={\mathbb P}(|X_0|>t)=\nu(|g|>t)$. From Subsection 2.2, Proposition \ref{prop:alphaW1}
and Item 3 of Proposition \ref{prop:suff}, we infer that $\sqrt n W_1(\tilde \mu_n, \mu)$ converges in distribution to the random variable
$\int |G(t)|dt$, where $G$ is a Gaussian random variable in ${\mathbb L}^1(dt)$ as soon as 
\beq\label{condH1}
\int_0^\infty (H(t))^{\frac{1-2 \gamma}{2(1-\gamma)}} dt < \infty \, .
\eeq
As a consequence:
\begin{enumerate}
\item If $g$ is positive and non increasing on (0, 1), with
\[
g(x) \leq \frac{C}{x^{(1-2\gamma)/2}|\ln(x)|^{b}}
\quad \text{near 0, for some $C>0$ and $b>1$,}
\]
then
\eref{condH1} holds.
\item If $g$ is positive and non decreasing on (0, 1), with
\[
g(x) \leq \frac{C}{(1-x)^{(1-2\gamma)/(2-2\gamma)}|\ln(1-x)|^{b}}\quad
\text{near 1, for some  $C>0$ and $b>1$,}
\]
then
\eref{condH1} holds.
\end{enumerate}
Recall from \eref{dual} that 
$W_1(\tilde \mu_n, \mu)= \sup_{f \in \Lambda_1}  | \tilde \mu_n(f) - \mu(f)|$, so
that the condition \eref{condH1} allows to control the supremum of $\sqrt n (\mu_n(f)-\mu(f))$ 
over the class $\Lambda_1$. Now if we only want a central limit theorem for 
$\sqrt n (\mu_n(f)-\mu(f))$ where $f$ is an element of $\Lambda_1$, then it follows from
\cite{DGM} that the condition 
\beq\label{condH2}
\int_0^\infty t(H(t))^{\frac{1-2 \gamma}{1-\gamma}} dt < \infty \, .
\eeq
is sufficient. For the two simple examples above, this would give the constraint $b>1/2$ instead
of $b>1$.

\subsection{Upper bounds for moments of order 1 and 2}
\label{M1/2}
In this section, we give some upper bounds for the 
quantities ${\mathbb E}(W_1(\mu_n, \mu))$ and 
$\|W_1(\mu_n, \mu)\|_2$ in terms of the 
coeffcients $\alpha_{1, \bf X}(k)$ and of the tail 
function $H$. 
For any $t \geq 0$, let
\begin{equation}\label{def:Sn}
S_{\alpha,n}(t)= \sum_{k=0}^n \min\left \{\alpha_{1, \bf X}(k), H(t)\right \} \, .
\end{equation}

\begin{prop}\label{M1and2} 
The following upper bounds hold:
\begin{equation}
{\mathbb E}(W_1(\mu_n, \mu))
\leq 4 \int_0^\infty \sqrt{\min \Big \{ \big(H(t)\big)^2, \frac{S_{\alpha, n}(t)}{n} \Big \}} \,  dt \, ,
\end{equation}
and
\begin{equation}\label{M2}
\|W_1(\mu_n, \mu)\|_2 \leq \frac{2 \sqrt 2}{\sqrt n} \int_0^\infty \sqrt{
S_{\alpha, n}(t) }
\,  dt \, .
\end{equation}
\end{prop}

\begin{rem}
As will be clear from the proof, one can also get
some upper bounds involving the quantity
$B(t)=F(t)(1-F(t))$ instead of $H(t)$. For instance, we 
can obtain an extension of the upper bound given in 
Theorem  3.5 of \cite{BL} to $\alpha$-dependent sequences. 
We have chosen to express the upper bounds in 
terms of the function $H$, because they are easier to
compute in the $\alpha$-dependent case (see 
Remark \ref{appl} below). 

The proof of Proposition \ref{M1and2} is based on the following
elementary inequality applied to $p=1$ and $p=2$:
$$
\text{For any  $p\geq 1$}, \quad 
\left \| \int | F_n(t)-F(t)| \ dt \right \|_p \leq 
\int \| F_n(t)-F(t)\|_p\  dt \, .
$$
One could also start from this inequality in the case where 
$p \in (1, 2)$ (resp. $p>2$) by applying a von Bahr-Esseen bound (resp. 
a Rosenthal bound) to $\| F_n(t)-F(t)\|_p$. 
However,  this would give less satisfactory bounds than in 
Subsections \ref{Sec:VBE} and \ref{Sec:Ros}, even in the i.i.d. case. For instance,  in the i.i.d. case and  $p \in (1,2)$, 
this would give 
\beq \label{badVBE}
\|W_1(\mu_n, \mu))\|^p_p \ll \frac{1}{n^{p-1}} \left (\int 
\left(H(t)\right)^{1/p} dt \right)^p \, .
\eeq
Note that the condition $\int 
(H(t))^{1/p} dt < \infty$  is more restrictive than $\|X_0\|_p< \infty$. Hence  the upper bound \eref{mdep} of Subsection 
\ref{Sec:VBE}
 is always better than
\eref{badVBE}.
\end{rem}

\begin{rem}
Starting from  Inequality \eref{Ebr}  and following the proof 
of Proposition \ref{M1and2} we obtain the upper bounds
\begin{equation}
{\mathbb E}(W_r^r(\mu_n, \mu))
\leq 4 \int_0^\infty t^{r-1}\sqrt{\min \Big \{ \big(H(t)\big)^2, \frac{S_{\alpha, n}(t)}{n} \Big \}} \,  dt \, ,
\end{equation}
and
\begin{equation}
\|W_r^r(\mu_n, \mu)\|_2 \leq \frac{2 \sqrt 2}{\sqrt n} \int_0^\infty t^{r-1}\sqrt{S_{\alpha, n}(t)}
\,  dt \, .
\end{equation}
\end{rem}

\begin{rem}\label{appl}
As a consequence of Proposition \ref{M1and2},  the following upper bounds hold:
\begin{enumerate}
\item
If \eref{condalpha} holds, 
then $\|W_1(\mu_n, \mu)\|_2 \ll n^{-1/2}$.
\item If $\alpha(k)=O(k^{-a})$ for some $a>1$, then
\begin{equation}\label{W1s}
{\mathbb E}(W_1(\mu_n, \mu))
\ll  \left( \int_0^{n^{-\frac{a}{a+1}}} Q(u) du + \frac{1}{\sqrt n} \int_{n^{-\frac{a}{a+1}}}^1
\frac{Q(u)}{u^{\frac{a+1}{2a}}} \, du \right)\, ,
\end{equation}
and
\begin{equation}\label{W2s}
\|W_1(\mu_n, \mu)\|_2 \ll \left( \int_0^{n^{-a}} \frac{Q(u)}{\sqrt u} du + \frac{1}{\sqrt n} \int_{n^{-a}}^1
\frac{Q(u)}{u^{\frac{a+1}{2a}}} \, du \right)\, .
\end{equation}
\item
If $\alpha(k)=O(a^k)$ for some $a<1$, then
\begin{equation*}
{\mathbb E}(W_1(\mu_n, \mu))
\ll \left( \int_0^{\frac{\ln (n)}{n}} Q(u) du + \frac{1}{\sqrt n} \int_{\frac{\ln(n)}{n}}^1
\frac{Q(u)|\ln (u)|}{\sqrt u} \, du \right)\, ,
\end{equation*}
and
\begin{equation*}
\|W_1(\mu_n, \mu)\|_2 \ll \left( \int_0^{e^{-n}} \frac{Q(u)}{\sqrt u} du + \frac{1}{\sqrt n} \int_{e^{-n}}^1
\frac{Q(u)|\ln (u)|}{\sqrt u} \, du \right)\, .
\end{equation*}
\item Assume that the $\alpha_k$'s converge to zero, but are not summable, and let
\begin{equation*}
u_n= \frac 1 n \sum_{k=1}^n \alpha_k \, .
\end{equation*}
Then
\begin{equation}\label{W1ns}
{\mathbb E}(W_1(\mu_n, \mu))
\ll \int_0^{\sqrt{u_n}} Q(u) \,  du \, ,
\end{equation}
and
\begin{equation}\label{W2ns}
\|W_1(\mu_n, \mu)\|_2  \ll \int_0^{u_n} \frac{Q(u)}{\sqrt u} \,  du \, .
\end{equation}
\end{enumerate}
\end{rem}

\begin{rem}
In the $m$-dependent case, the inequality \eref{W1s} holds with $a=\infty$, that is
$$
{\mathbb E}(W_1(\mu_n, \mu))
\ll \left( \int_0^{n^{-1}} Q(u) du + \frac{1}{\sqrt n} \int_{n^{-1}}^1
\frac{Q(u)}{\sqrt{u}} \, du \right)\, .
$$
In particular, if $H(t)=O(t^{-1} (\ln (t))^{-a})$ for some $a>1$ (which implies that
${\mathbb E}(|X_0|)< \infty$), then
$Q(u) =O(u^{-1}|\ln(u)|^{-a})$,
and consequently
$$
{\mathbb E}(W_1(\mu_n, \mu)) \ll \frac{1}{(\ln (n))^{a-1}} \, .
$$
\end{rem}

\noindent {\bf Example (continued).}  We continue the example of Subsection \ref{CLT}. 
\begin{enumerate}
\item If $g$ is positive and non increasing on (0, 1), with
\[
g(x) \leq \frac{C}{x^{b}}
\quad \text{near 0, for some $C>0$ and  $b\in [0,  1-\gamma)$,}
\]
then $Q(u)\leq D u^{-b/(1-\gamma)}$ for some $D>0$.  Applying  \eref{W1s}-\eref{W2s} and  \eref{W1ns}-\eref{W2ns}, the following upper bounds hold.

For  $\gamma \in (0, 1/2)$, 
\beq 
{\mathbb E}(W_1(\tilde \mu_n, \mu))
\ll
\begin{cases}
   n^{-1/2}  \quad \quad \quad \ \text{if $b<(1-2\gamma)/2$}\\
   n^{-1/2}  \ln (n)  \quad \text{if $b=(1-2\gamma)/2$} \\
   n^{b+\gamma -1}  \quad \quad \ \, \, \, 
  \text{if $ b > (1-2\gamma)/2$,}
  \end{cases}
  \nonumber
\eeq 
and 
\beq
\|W_1(\tilde \mu_n, \mu)\|_2
\ll
\begin{cases}
  n^{-1/2} \quad \quad \quad \ \, \text{ if $b<(1-2\gamma)/2$}\\
  n^{-1/2}  \ln (n) \  \, \quad  \text{if $b=(1-2\gamma)/2$} \\
  n^{(2b+\gamma -1)/2\gamma} \quad  \text{if $(1-2\gamma)/2< b < (1-\gamma)/2$.}
  \end{cases}
  \nonumber
\eeq

For $\gamma=1/2$, 
$$
{\mathbb E}(W_1(\tilde \mu_n, \mu)) \ll \left(\frac{\ln(n)}{n}\right)^{\frac{1-2b}2},  \ \text{and}  \
\|W_1(\tilde \mu_n, \mu)\|_2 
\ll \left(\frac{\ln(n)}{n}\right)^{\frac{1-4b}2}   \text{if $b<1/4$.}
$$

For $\gamma \in (1/2,1)$, 
$$
{\mathbb E}(W_1(\tilde \mu_n, \mu)) \ll  n^{\frac{b+\gamma-1}{2\gamma}},  \ \text{and}  \
\|W_1(\tilde \mu_n, \mu)\|_2 
 \leq C n^{\frac{2b+\gamma -1}{2\gamma}}  \ \text{if $b<(1-\gamma)/2$.}
$$
\item If $g$ is positive and non decreasing on (0, 1), with
\[
g(x) \leq \frac{C}{(1-x)^{b}}
\quad \text{near 1, for some $C>0$ and $b \in [0, 1)$,}
\]
then $Q(u)\leq D u^{-b}$ for some $D>0$.  Applying  \eref{W1s}-\eref{W2s} and  \eref{W1ns}-\eref{W2ns}, the following upper bounds hold.

For  $\gamma \in (0, 1/2)$, 
\beq 
{\mathbb E}(W_1(\tilde \mu_n, \mu))
\ll
\begin{cases}
  n^{-1/2}  \quad \quad \quad \ \text{if $b<(1-2\gamma)/2(1-\gamma)$}\\
    n^{-1/2}  \ln (n)  \quad \text{if $b=(1-2\gamma)/2(1-\gamma)$} \\
  n^{(\gamma-1)(1-b)}  \quad \,
  \text{if $ b > (1-2\gamma)/2(1-\gamma)$,}
  \end{cases}
  \nonumber
\eeq 
and 
\beq
\|W_1(\tilde \mu_n, \mu)\|_2
\ll
\begin{cases}
  n^{-1/2} \quad \quad \quad \quad \ \, \text{ if $b<(1-2\gamma)/2(1-\gamma)$}\\
   n^{-1/2}  \ln (n) \  \, \quad \quad  \text{if $b=(1-2\gamma)/2(1-\gamma)$} \\
  n^{(\gamma-1)(1-2b)/2\gamma} \quad  \text{if $(1-2\gamma)/2(1-\gamma)< b < 1/2$.}
  \end{cases}
  \nonumber
\eeq
For $\gamma=1/2$, 
$$
{\mathbb E}(W_1(\tilde \mu_n, \mu)) \ll \left(\frac{\ln(n)}{n}\right)^{\frac{1-b}2},  \ \text{and}  \
\|W_1(\tilde \mu_n, \mu)\|_2 
\ll \left(\frac{\ln(n)}{n}\right)^{\frac{1-2b}2}   \text{if $b<1/2$.}
$$
For $\gamma \in (1/2,1)$, 
$$
{\mathbb E}(W_1( \tilde \mu_n, \mu)) \ll  n^{\frac{(\gamma-1)(1-b)}{2\gamma}},  \ \text{and}  \
\|W_1( \tilde\mu_n, \mu)\|_2 
\ll n^{\frac{(\gamma -1)(1-2b)}{2\gamma}}  
\ \text{if $b<1/2$.}
$$
\end{enumerate}

\noindent {\bf Proof of Proposition \ref{M1and2}.} 
Starting from \eref{empL1},  we immediately see that
\begin{equation}\label{b0}
{\mathbb E}(W_1(\mu_n, \mu)) \leq \int \|F_n(t)-F(t)\|_1 \ dt  \quad
\text{and} 
\quad 
\|W_1(\mu_n, \mu)\|_2 \leq \int  \|F_n(t)-F(t)\|_2 \ dt \, .
\end{equation}
Let $B(t)=F(t)(1-F(t))$, and note first that
\begin{equation}\label{b1}
\|F_n(t)-F(t)\|_1 \leq \|{\bf 1}_{X_0\leq t}-F(t)\|_1 = 2 B(t) \, .
\end{equation}
On another hand
\begin{equation}\label{b2}
\|F_n(t)-F(t)\|_1^2 \leq \|F_n(t)-F(t)\|_2^2 \leq 
\frac 1 n {\mathrm {Var}}({\bf 1}_{X_0\leq t}) +
\frac 2 n \sum_{k=1}^n \left | 
{\mathrm {Cov}}({\bf 1}_{X_0\leq t}, {\bf 1}_{X_k\leq t})
 \right |  \,  . 
\end{equation}
Now, the two following upper bounds hold:
\begin{align}
\label{b3}
\left |{\mathrm {Cov}}({\bf 1}_{X_0\leq t}, {\bf 1}_{X_k\leq t})\right | &\leq 
  \|{\mathbb E}({\bf 1}_{X_k\leq t}|{\mathcal F}_0)-F(t)\|_1 \leq \alpha_{1, \bf X}(k)  \, , \\
  \label{b4}
\left |{\mathrm {Cov}}({\bf 1}_{X_0\leq t}, {\bf 1}_{X_k\leq t})
\right |   
&\leq  {\mathrm {Var}}({\bf 1}_{X_0\leq t})=B(t) \, .
\end{align}
From \eref{b1}, \eref{b2}, \eref{b3} and \eref{b4} 
it follows that 
$$
\|F_n(t)-F(t)\|_1 \leq 2 \sqrt{\min \Big \{ \big(B(t)\big)^2,
\frac 1 n  \sum_{k=0}^n \min\left \{\alpha_{1, \bf X}(k), B(t)\right \} \Big \}}
$$
and 
$$
\|F_n(t)-F(t)\|_2 \leq \sqrt{ \frac {2} {n} 
  \sum_{k=0}^n \min\left \{\alpha_{1, \bf X}(k), B(t)\right \} } \, .
$$
These two upper bounds combined with \eref{b0} imply that 
\begin{multline*}
{\mathbb E}(W_1(\mu_n, \mu)) \leq 2 \int \sqrt{\min \Big \{ \big(B(t)\big)^2,
\frac 1 n  \sum_{k=0}^n \min\left \{\alpha_{1, \bf X}(k), B(t)\right \} \Big \}} \ dt \\
\leq 4 
\int_0^\infty \sqrt{\min \Big \{ \big(H(t)\big)^2, \frac{S_{\alpha, n}(t)}{n} \Big \}} \  dt
\end{multline*}
and 
$$
\|W_1(\mu_n, \mu)\|_2  \leq \sqrt {\frac 2  n} \int \sqrt{
  \sum_{k=0}^n \min\left \{\alpha_{1, \bf X}(k), B(t)\right \} } \ dt \leq 
\frac{2 \sqrt 2 }{\sqrt n} \int_0^\infty \sqrt{
S_{\alpha, n}(t) }
\  dt \, ,  
$$
which are the desired inequalities. 
\subsection{A von Bahr-Esseen type inequality}\label{Sec:VBE}
In this section, we give some upper bounds for the 
quantity $\|W_1(\mu_n, \mu)\|_p$ when $p \in (1,2)$ 
 in terms of the 
coefficients $\alpha_{1, \bf X}(k)$ and of the quantile function $Q$. 
For $u\in (0,1)$, let
\begin{equation} \label{alphamoins1}
\alpha_{1,{\bf X}}^{-1}(u)=\sum_{k=0}^\infty
 {\bf 1}_{u \leq \alpha_{1,{\bf X}}(k)} \, .
\end{equation}

\begin{prop}\label{propvBE}
For $p  \in (1,2)$, the following inequality holds
\beq \label{vBEineqW1}
\|  W_1 (\mu_n , \mu) \Vert_p^p \ll 
\frac 1 {n^{p-1}}  \int_0^1 (\alpha_{1,{\bf X}}^{-1} (u) \wedge n)^{p-1}   Q^{p}(u) du \, .
\eeq
\end{prop}
Note that Inequality \eref{vBEineqW1} writes also
$$
\|  W_1 (\mu_n , \mu) \Vert_p^p \ll 
\frac 1 {n^{p-1}}  \sum_{k=0}^n \frac 1 {(k+1)^{2-p}}\int_0^{\alpha_{1, \bf X}(k)} Q^{p}(u) du \, .
$$
\begin{rem} Let $r \geq 1$ and $p  \in (1,2)$. Starting again from \eref{Ebr}
and following the proof of Proposition \ref{propvBE}, we obtain the 
upper bound
\beq \label{vBEineqWr}
\Vert  W_r^r (\mu_n , \mu) \Vert_p^p \ll 
 \frac 1 {n^{p-1}} \int_0^1 (\alpha_{1,{\bf X}}^{-1} (u) \wedge n)^{p-1}   Q^{pr}(u) du \, .
\eeq
\end{rem}

\begin{rem}
In the $m$-dependent case, Inequality \eref{vBEineqWr}
becomes
\beq \label{mdep}
\| W_1 (\mu_n , \mu) \|_p^p \ll 
\frac 1 {n^{p-1}}  \|X_0\|_p^p \, .
\eeq
This inequality seems to be new even in the i.i.d. case. 
It is noteworthy  that the upper bound \eref{mdep}
is the same as the moment bound of order $p$ for partial sums
of i.i.d. random variables, which can be deduced from the
 classical inequality of von Bahr and Esseen 
\cite{BE}.
\end{rem}

\noindent {\bf Example (continued).}  We continue the example of Subsection \ref{CLT}. 
\begin{enumerate}
\item Let $p \in (0,1)$, and let  $g$ be positive and non increasing on (0, 1), with
\[
g(x) \leq \frac{C}{x^{b}}
\quad \text{near 0, for some $C>0$ and $b\in [0, (1-\gamma)/p)$.}
\]
 Applying Proposition  \ref{propvBE}, the following upper bounds hold.

For  $\gamma \in (0, 1/p)$, 
\beq 
\|W_1(\tilde \mu_n, \mu))\|_p
\ll
\begin{cases}
   n^{(1-p)/p} \quad \quad \quad   \quad \ \ \,  \text{if $b<(1-p\gamma)/p$}\\
 (n^{(1-p)}  \ln (n))^{1/p}    \ \quad \text{if $b=(1-p\gamma)/p$} \\
   n^{(pb + \gamma -1)/p \gamma} 
 \quad \quad  \quad  
  \text{if $ b > (1-p\gamma)/p$.}
  \end{cases}
  \nonumber
\eeq 
Moreover, if $b=(1-p\gamma)/p$, Proposition \ref{propvBEine} below gives the upper bound
\begin{equation}\label{better}
{\mathbb P} \left ( W_1 (\mu_n , \mu)  \geq   x \right )
\ll \frac{1}{n^{p-1} x^p} \, .
\end{equation}

For $\gamma \in [1/p, 1)$, 
$
\|W_1(\tilde \mu_n, \mu))\|_p 
\ll  n^{(pb + \gamma -1)/p \gamma}
$.
\item Let $p \in (0,1)$, and let  $g$ be positive and non decreasing on (0, 1), with
\[
g(x) \leq \frac{C}{(1-x)^{b}}
\quad \text{near 1, for some $C>0$ and $b\in [0, 1/p)$.}
\]
 Applying Proposition  \ref{propvBE}, the following upper bounds hold.

For  $\gamma \in (0, 1/p)$, 
\beq 
\|W_1(\tilde \mu_n, \mu))\|_p
\ll
\begin{cases}
   n^{(1-p)/p} \quad \quad \quad   \quad \ \ \,  \text{if $b<(1-p\gamma)/(p(1-\gamma))$}\\
 (n^{(1-p)}  \ln (n))^{1/p}   \ \quad  \text{if $b=(1-p\gamma)/(p(1-\gamma))$} \\
   n^{(\gamma-1)(1-pb)/p \gamma}  \quad  \ \
  \text{if $ b > (1-p\gamma)/(p(1-\gamma))$.}
  \end{cases}
  \nonumber
\eeq 

Moreover, if $b=(1-p\gamma)/(p(1-\gamma))$, Proposition \ref{propvBEine} below gives the upper bound
\eref{better}.

\medskip

For $\gamma \in [1/p, 1)$, 
$
\|W_1(\tilde \mu_n, \mu))\|_p \ll
 n^{(\gamma-1)(1-pb)/p \gamma} 
$.

\end{enumerate}

\begin{rem}
The upper bound \eref{better} is in accordance with a result by Gou\"ezel \cite{Gou}.  He  proved that, if 
$g$ is exactly of the form
$g(x)=x^{-(1-p\gamma)/p}$ and $\theta$ is the LSV map defined by \eref{LSVmap}, then for any 
positive real $x$,
$$
\lim_{n \rightarrow \infty} \nu
\left(\frac{1}{n^{1/p}}\left | \sum_{k=1}^n  \left (g \circ \theta^k - \nu(g) \right ) \right| > x \right)=
{\mathbb P}(|Z_p|>x) \, ,
$$
where $Z_p$ is a $p$-stable random variable such that $\lim_{x \rightarrow \infty}
x^p{\mathbb P}(|Z_p|>x)=c>0$.
\end{rem}

\medskip

\noindent {\bf Proof of Proposition \ref{propvBE}}. 
For any $n \in {\mathbb N}$, let us introduce the following notations:
  \begin{equation*}
  R_n(u)=(\min\{q\in {\mathbb N}^* \tq \alpha_{1,{\bf X}}(q)\leq u\} \wedge n) Q(u)
  \text{ and }
R_n^{-1}(x)=\inf\{u\in [0,1] \tq R_n(u)\leq x\}\,.
  \end{equation*}
The proof is based on the following proposition:
\begin{prop}\label{propvBEine}
For any positive integer $n$, any $x>0$, and any  $\eta \in [1, 2[$, the following inequality holds:
\beq \label{vBEineqW1bis}
{\mathbb P} \left ( n W_1 (\mu_n , \mu)  \geq  6 x \right ) \leq c_1 \frac{n}{x}  \int_0^{R_n^{-1}(x)}Q(u) du + c_2 \frac{n}{x^{\eta}}  \int_{R_n^{-1}(x)}^1 R^{\eta-1}_n(u)Q(u) du  \, ,
\eeq
where $c_1 = 36$ and $c_2 = 64(2-\eta)^{-1} $.
\end{prop}
Before proving the proposition above, let us see how it entails Proposition \ref{propvBE}. We have 
\[
\Vert n W_1 (\mu_n , \mu) \Vert_p^p  = 6^p p \int_0^{\infty} x^{p-1} {\mathbb P} \big ( n W_1 (\mu_n , \mu)  \geq  6 x \big ) dx \, .
\]
Therefore applying Inequality \eqref{vBEineqW1bis} with 
$\eta \in (p, 2)$ and using the fact that
\[
u < R_n^{-1}  (x) \iff x < R_n(u) \, ,
\]
we get 
\begin{multline*}
\Vert n W_1 (\mu_n , \mu) \Vert_p^p  \leq 6^p p \, n c_1 \int_0^1 Q(u)  \int_0^{\infty} x^{p-2} \I_{x < R_n(u)}\  dx du \\
+ 6^p p\,  n c_2 \int_0^1 R_n^{\eta-1}(u)Q(u)  \int_0^{\infty} x^{p-1 - \eta } \I_{x \geq  R_n(u)} \  dx du \, ,
\end{multline*}
which gives the desired result since $1 < p < \eta < 2$. 
Hence it remains to prove Proposition \ref{propvBEine}.
\medskip

\noindent {\bf Proof of Proposition \ref{propvBEine}}. Let 
\beq \label{defvM}v = R_n^{-1}(x) \, , \, M=Q(v) 
\eeq
and set $ g_M(y) = (y \wedge M) \vee (-M) $.  For any integer $i$, let
\beq \label{detrun}
X_i' = g_M(X_i )\,  \mbox{ and } \, X_i'' = X_i -X_i' \, .
\eeq
Starting from \eref{dual}, we first notice that
\begin{align*}
n W_1 (\mu_n , \mu) &= \sup_{f \in \Lambda_1} \sum_{i=1}^n \left ( f(X_i) - \E(f(X_i)) \right )\\
& \leq  
\sup_{f \in \Lambda_1} \sum_{i=1}^n 
\left ( f(X'_i) - \E(f(X'_i)) \right )
+
\sup_{f \in \Lambda_1} \sum_{i=1}^n 
\left ( f(X_i) - f(X_i') -  \E(f(X_i)-f(X_i')) \right )\, .
\end{align*}
Therefore
\beq \label{decW1Wprime}
n W_1 (\mu_n , \mu)  \leq \sup_{f \in \Lambda_1} \sum_{i=1}^n \left(f(X'_i) - \E(f(X'_i)) \right)
+ \sum_{i=1}^n ( |X_i''| + \E( |X_i''|) \, .
\eeq
Let now 
\beq \label{defofq}
q = \min\{k\in {\mathbb N}^* \tq \alpha_{1,{\bf
X}}(k)\leq v\} \wedge n \, . 
\eeq
Since $R_n$ is right
continuous, we have $R_n(R^{-1}_n(w))\leq w$ for any $w$, hence
\beq \label{restq1}
qM = R_n(v) = R_n (R_n^{-1}(x)) \leq x \, .
\eeq
Assume first that $q=n$. Bounding  $f(X'_i) - \E(f(X'_i)) $ by $2M$ in \eqref{decW1Wprime}, we obtain
\begin{equation}\label{vB1}
n W_1 (\mu_n , \mu)   \leq 2qM + \sum_{k=1}^n ( |X''_k| + \E( |X_k''|) )\, .
\end{equation}
Taking into account \eqref{restq1} this gives
  \begin{equation*}
 {\mathbb P} \big ( n W_1 (\mu_n , \mu)  \geq 6x \big )  
  \leq \frac{1}{2x} \sum_{k=1}^n  \E ( | X_k''|) \, .
  \end{equation*}
Writing $\varphi_M(x)=(|x|-M)_+$, we have
\[
\sum_{k=1}^n \E ( | X_k''|) \leq  \sum_{k=1}^n \E (\varphi_{M} ( X_k)) \, .
\]
But $Q_{\varphi_M (X_k)} \leq Q_{| X_k |}\I_{[0,v]} \leq Q\I_{[0,v]} $. Consequently
   \begin{equation} \label{dec13FN}
  \sum_{k=1}^n  \E ( | X_k''|)
  \leq  n \int_0^{R^{-1}_n(x)} Q(u) du\, .
  \end{equation}
From \eref{vB1} and \eref{dec13FN}, we infer that  \begin{equation} \label{dec13FNdef}
  {\mathbb P} \big ( n W_1 (\mu_n , \mu)  \geq 6x \big )  
  \leq \frac{n}{2x}  \int_0^{R^{-1}_n(x)} Q(u) du  \, ,
  \end{equation}
which then proves the proposition in case where $q=n$.

From now on, we assume that $q<n$. Therefore $q =  \min\{k\in {\mathbb N}^* \tq \alpha_{1,{\bf
X}}(k)\leq v\} $ and then $\alpha_{1,{\bf
X}} (q)\leq v$. Starting from 
\eqref{decW1Wprime}, we first notice that
\[ 
 {\mathbb P} \left  ( n W_1 (\mu_n , \mu)  \geq 6x \right )  
  \leq   {\mathbb P} \left ( \sup_{f \in \Lambda_1} \sum_{i=1}^n\left ( f(X'_i) - \E(f(X'_i))\right ) \geq 5 x \right )  
  + \frac{2}{x} \sum_{k=1}^n  \E ( | X_k''|) \, .
\]
Therefore taking into account \eqref{dec13FN},
 \begin{equation} \label{dec14FN}
   {\mathbb P} \left ( n W_1 (\mu_n , \mu)  \geq 6x 
   \right )  
  \leq   {\mathbb P} \left ( \sup_{f \in \Lambda_1} \sum_{i=1}^n \left ( f(X'_i) - \E(f(X'_i)) \right )\geq 5 x \right )  +  \frac{2n}{x} \int_0^{R^{-1}_n(x)} Q(u) du\, .
  \end{equation}
To control the first term on the right-hand side, we first notice that 
\begin{multline*}
 \sup_{f \in \Lambda_1} \sum_{i=1}^n 
 \left (f(X'_i) - \E(f(X'_i))\right ) \leq   \sup_{f \in \Lambda_1} \sum_{i=1}^{[n/q]q} \left ( f(X_i') - {\mathbb E} (f(X_i'))\right ) + 2 (n-[n/q]q) M \\
\leq   \int_{\mathbb R} \left \vert \sum_{i=1}^{[n/q]q} \left (\I_{X_i' \leq t} - {\mathbb E} \left (\I_{X_i' \leq t} 
\right) \right)\right \vert dt + 2 q M .
\end{multline*}
Using \eqref{restq1}, it follows that 
\begin{equation*}
{\mathbb P} \left ( \sup_{f \in \Lambda_1} \sum_{i=1}^n f(X'_i) - \E(f(X'_i)) \geq 5 x \right )   
\leq
{\mathbb P} \left ( \int_{\mathbb R} \left \vert \sum_{i=1}^{[n/q]q}\left(\I_{X_i' \leq t} - {\mathbb E} \left (\I_{X_i' \leq t} \right ) \right)\right \vert dt  \geq 3 x \right )  \, .
\end{equation*}
For any integer $i$, define
\[ U_i(t) =  \sum_{k=(i-1)q+1}^{iq} \left(\I_{X_k' \leq t} - {\mathbb E}\left  (\I_{X_k' \leq t} \right ) \right ) \, .\]
Consider now the $\sigma$-algebras ${\cal G}_{i}= \F_{iq}$ and define the  variables $\tilde
U_i(t)$ as follows:   $\tilde U_{2i -1} (t) =  U_{2i -1}(t)-
\E (U_{2i - 1 } (t) | {\cal G}_{2(i-1) -1 })$  and  $ \tilde U_{2i}(t)=U_{2i}(t)- \E (U_{2i} (t)| {\cal
G}_{2(i-1)}) $. Substituting $\tilde U_i(t)$ to $U_i(t)$, we
obtain the inequality
\begin{multline} \label{decavecmart}
\left \vert \sum_{i=1}^{[n/q]q} \left(\I_{X_i' \leq t} - {\mathbb E} \left (\I_{X_i' \leq t} \right ) \right ) \right \vert = \left  \vert \sum_{i=1}^{[n/q]} U_i(t) \right \vert \\
 \leq  \max_{2 \leq 2j \leq
 [n/q]
 } \left \vert \sum_{i=1}^j \tilde U_{2i} (t) \right \vert \!+ \!\max_{1 \leq 2j-1 \leq
 [n/q]
 } \left \vert \sum_{i=1}^j \tilde U_{2i -1} (t) \right \vert
 \!+\!
 \sum_{i=1}^{ [n/q]
 } \vert  U_{i} (t) -\tilde U_{i} (t) \vert \, .
\end{multline}
Therefore 
\begin{equation} \label{decavecmartproba}
{\mathbb P} \left ( \sup_{f \in \Lambda_1} \sum_{i=1}^n \left ( f(X'_i) - \E(f(X'_i)) \right ) \geq 5 x \right )   
\leq  I_1(n) + I_2(n) + I_3(n) \, ,
\end{equation}
where
\begin{align*}
I_1(n)&={\mathbb P} \left ( \int_{\mathbb R}  \sum_{i=1}^{ [n/q]
 } \vert  U_{i} (t) -\tilde U_{i} (t) \vert  \ dt \geq  x \right )  \\
I_2(n) & = {\mathbb P} \left ( \int_{\mathbb R} \max_{2 \leq 2j \leq
 [n/q]
 } \left \vert \sum_{i=1}^j \tilde U_{2i} (t) \right \vert \ dt  \geq  x \right )  \\
I_3(n) & =   {\mathbb P} \left ( \int_{\mathbb R} \max_{1 \leq 2j-1 \leq
 [n/q]
 } \left \vert \sum_{i=1}^j \tilde U_{2i-1} (t)\right \vert \ dt  \geq  x \right ) \, .
\end{align*}
Using Markov's inequality and stationarity, we get
\begin{equation*} I_1(n) \leq \frac{n}{x} \int_{\mathbb R} {\mathbb E} \left\vert  {\mathbb E}\left (\I_{X_1' \leq t}  | \F_{-q} \right)  - \E \left (\I_{X_1' \leq t} \right ) \right \vert
\ dt   = \frac{n}{x} \int_{-M}^M {\mathbb E} \left \vert  {\mathbb E}\left (\I_{X_1' \leq t}  | \F_{-q} \right )  - \E \left (\I_{X_1' \leq t} \right ) \right \vert \ dt    \, .
\end{equation*}
But, 
$$
\sup_{t \in {\mathbb R} }\left \Vert  {\mathbb E}
\left (\I_{X_1' \leq t}  | \F_{-q} \right )  - \E 
\left (\I_{X_1' \leq t} \right ) \right \Vert_1 = \sup_{t \in {\mathbb R} } 
\left \Vert  {\mathbb E}\left (\I_{g_M(X_1) \leq  t}  | \F_{-q} 
\right )  - \E \left (\I_{g(X_1) \leq  t} \right ) \right \Vert_1
\leq \alpha_{1,{\bf
X}} (q+1) \, ,
$$
where the inequality comes from the fact that $g_M$ is a nondecreasing function. Therefore, 
\begin{equation} \label{I1n}
I_1(n) 
\leq \frac{2n}{x}   Q(v) \alpha_{1,{\bf
X}} (q+1)  \leq \frac{2n}{x}  v  Q(v) \leq \frac{2n}{x}  \int_0^v Q(u) du   \, .
\end{equation}
We handle now the term $I_2(n)$ in the decomposition \eqref{decavecmartproba}. Using again Markov's inequality, we get 
\[
I_2(n)  \leq \frac 1 {x^2} \left ( \int_{-Q(v)}^{Q(v)}   
\left \| \max_{2 \leq 2j \leq
 [n/q] } \left | \sum_{i=1}^j \tilde U_{2i} (t) \right |
\right  \|_2 dt   \right )^2  \, .
\]
By Doob's maximal inequality,  
\[ \left \Vert \max_{2 \leq 2j \leq
 [n/q]
 } \left \vert \sum_{i=1}^j \tilde U_{2i} (t) \right \vert \right \Vert^2_2
\leq 2 \sum_{i=1}^{\frac12  [n/q] } \Vert  \tilde U_{2i} (t) \Vert_2^2  \leq 2 \sum_{i=1}^{\frac12  [n/q] } \Vert   U_{2i} (t) \Vert_2^2 \, .
\]
Now 
\begin{multline*}   
 2 \sum_{i=1}^{\frac12  [n/q] } \Vert   U_{2i} (t) \Vert_2^2 \leq \frac{n}{q} 
 \left\Vert   \sum_{k=1}^{q} 
 \left (\I_{X_k' \leq t}- {\mathbb E} \left ( \I_{X_k' \leq t} 
 \right ) \right )   \right \Vert_2^2 \\
 \leq 2n  \sum_{k=0}^{q-1} \left | \E 
 \left ( \left (\I_{g_M(X_0) \leq t}- {\mathbb E} \left ( \I_{g_M(X_0)  \leq t} \right )\right  ) \left  (\I_{g_M(X_k) \leq t}- {\mathbb E} \left ( \I_{g_M(X_k)  \leq t} \right )
 \right  ) \right  )  \right  | \, .
\end{multline*}
Note that since $g_M$ is a nondecreasing function, 
\begin{multline*}   
\sup_{t \in {\mathbb R}} \left | \E \left ( 
\left (\I_{g_M(X_0) \leq t}- {\mathbb E}\left ( \I_{g_M(X_0)  \leq t}\right  ) \right  ) \left  (\I_{g_M(X_k) \leq t}- {\mathbb E} 
\left ( \I_{g_M(X_k)  \leq t}\right  ) \right )\right  )  \right  | \\
\leq \sup_{t \in {\mathbb R}} \left | \E \left  ( 
\left (\I_{X_0 \leq  t}- {\mathbb E}\left  ( \I_{X_0   \leq t} \right )
\right  )  \left (\I_{X_k \leq t}- {\mathbb E} \left ( \I_{X_k  \leq t} \right ) \right ) \right ) 
 \right  | \, .
\end{multline*}
Moreover 
\begin{multline*}
 \sup_{t \in {\mathbb R}} \left | \E \left ( 
 \left (\I_{X_0> t}- {\mathbb E} \left ( \I_{X_0  > t} \right
 ) \right )  \left (\I_{X_k > t}- {\mathbb E}\left  ( \I_{X_k  > t} \right )\right ) \right  )  \right | \\ \leq \sup_{t \in {\mathbb R} } \left \Vert  {\mathbb E}
 \left (\I_{X_k> t}  | \F_{0} \right )  - \E \left (\I_{X_k > t}\right  ) \right \Vert_1 = \alpha_{1,{\bf
X}} (k)\, .
\end{multline*}
On an other hand, the following bound is also valid
\begin{multline*}
 \left | \E \left ( \left (\I_{g_M(X_0) \leq t}- {\mathbb E} 
 \left ( \I_{g_M(X_0)  \leq t} \right) \right ) 
 \left  (\I_{g_M(X_k) \leq t}- {\mathbb E}\left  ( \I_{g_M(X_k)  \leq t}\right  ) 
 \right ) \right )  \right | \\
 \leq {\mathrm {Var}}\left (\I_{g_M(X_0) \leq  t}\right ) \leq  \min \left 
 \{\E \left  (\I_{g_M(X_0) \leq  t} \right) , \E  
 \left (\I_{g_M(X_0) > t} \right ) 
 \right \} \, .
\end{multline*}
So, overall, we get
\begin{multline*}   
I_2(n)  \\
 \leq  \frac{2n} {x^{2}}  \left( \int_{0}^{Q(v)}    \left ( \sum_{k=0}^{q-1}   \alpha_{1,{\bf
X}} (k) \wedge {\mathbb P}  (g_M(X_0) > t )  \right )^{\frac 1 2 } 
 +     \left ( \sum_{k=0}^{q-1}  \alpha_{1,{\bf
X}} (k) \wedge {\mathbb P}  (- g_M(X_0) \geq  t )   \right )^{\frac 1 2 } dt    \right)^2 \\
\leq  \frac {2n} {x^{2}}  \left ( \int_{0}^{Q(v)}    \left ( \sum_{k=0}^{q-1}  \alpha_{1,{\bf
X}} (k) \wedge {\mathbb P}  (|X_0|> t )   \right )^{\frac 1 2 }  +   \left ( \sum_{k=0}^{q-1}   \alpha_{1,{\bf
X}} (k) \wedge {\mathbb P}  (|X_0| \geq  t )   \right )^{\frac 1 2 } dt    \right )^2  \, .
\end{multline*}
We then derive that
\begin{multline*}
I_2(n)  \leq  \frac{8 n} {x^{2} } \left ( \int_{0}^{Q(v)}    \left ( \sum_{k=0}^{q-1} \alpha_{1,{\bf
X}} (k) \wedge H(t) \right )^{1/2} dt   \right )^2 
\\ = 
\frac{8 n} {x^{2} } \left ( \int_{0}^{Q(v)}    \left ( \sum_{k=0}^{q-1} \int_0^{H(t)} \I_{u \leq \alpha_{1,{\bf
X}} (k) } du \right )^{1/2} dt   \right )^2
.
\end{multline*}
Using the fact that $\sum_{k=0}^{q-1} \I_{u \leq \alpha_{1,{\bf
X}} (k) } = \alpha^{-1}_{1,{\bf
X}} (u)  \wedge q $, we then get 
\begin{multline}   \label{I2n1}
I_2(n)  \leq  \frac{8 n} {x^{2} }  \left ( \int_{0}^{Q(v)}    \left (  v q +  \int_v^{H(t)} \left(  \alpha^{-1}_{1,{\bf
X}} (u) \wedge n\right ) du  \right )^{\frac 1 2 } dt   \right )^2 \\ \leq  
\frac{16 n} {x^{2} }  v q (Q(v))^2 + \frac{16 n} {x^{2} } \left ( \int_{0}^{Q(v)}    \left (  \int_v^{H(t)} \left(  \alpha^{-1}_{1,{\bf
X}} (u) \wedge n \right ) du  \right )^{\frac 1 2 } dt   \right )^2 \\
\leq \frac{16 n} {x } \int_0^v Q(u) du  + \frac{16 n} {x^{2} }  \left ( \int_{0}^{Q(v)}    \left (  \int_v^{H(t)} \left ( \alpha^{-1}_{1,{\bf
X}} (u) \wedge n \right) du  \right )^{1/2} dt   \right )^2 
\, .
\end{multline}
 where for the last inequality we have used \eqref{restq1} and the fact that  $v Q(v) \leq \int_0^v Q(u) du  $, since $Q$ is non increasing. To handle the last term on the right-hand side, we proceed as follows.  For any  $\eta$ in $[1, 2)$, we first note that  
\[
 \int_{v}^{H(t)} ( \alpha^{-1}_{1,{\bf
X}} (u) \wedge n ) du =  \int_{v}^{H(t)} \left ( \alpha^{-1}_{1,{\bf
X}} (u)  \wedge n  \right ) Q^{\eta} (u)Q^{-\eta} (u)  du \leq   \frac 1 {t^\eta}  \int_{v}^{H(t)}\left ( \alpha^{-1}_{1,{\bf
X}} (u)  \wedge n  \right )   Q^{\eta} (u)  du \, ,
\]
where the inequality comes from the fact that  $u < H(t) \iff t < Q(u)$, and then $u < H(t) $ implies that  $Q^{-\eta} (u) < t^{- \eta}$. Now, since $u >v$ 
implies  that  $\alpha^{-1}_{1,{\bf
X}} (u) \leq \alpha^{-1}_{1,{\bf
X}} (v) $, we get 
\[
\int_{v}^{H(t)} \left (  \alpha^{-1}_{1,{\bf
X}} (u) \wedge n \right ) du  \leq \frac 1 {t^\eta}\left (\alpha^{-1}_{1,{\bf
X}}  (v)  \wedge n \right )^{2 - \eta}  \int_{v}^{1} \left (\alpha^{-1}_{1,{\bf
X}}  (u)  \wedge n \right )^{ \eta -1} Q^{\eta} (u) du  \, .
\]
Therefore, since  $\eta/2 <1$,
\begin{multline*}
\left ( \int_{0}^{Q(v)}    \left (  \int_v^{H(t)}  \left (  \alpha^{-1}_{1,{\bf
X}} (u) \wedge n \right)  du  \right )^{1/2} dt   \right )^2  
\\ \leq   \left (\alpha^{-1}_{1,{\bf
X}}  (v) \wedge n \right )^{2 - \eta}  \left ( \int_0^{Q(v)} t^{-\eta	/2} dt \right)^2  \int_v^1  \left ( R_n^{-1}(u) \right )^{ \eta -1}  Q(u) du \\
\leq  \frac{2}{2-\eta} \left (\alpha^{-1}_{1,{\bf
X}}  (v) \wedge n \right )^{2 - \eta} Q^{2 - \eta} (v) \int_v^1  \left ( R_n^{-1}(u) \right)^{ \eta -1}  Q(u) du \, .
\end{multline*}
But, by \eqref{restq1}, $ \left (\alpha_{1,{\bf
X}} ^{-1} (v) \wedge n \right )^{2 - \gamma}  Q^{2- \gamma} (v)  = R_n^{2-\gamma}(v) \leq x^{2 - \gamma}$. Therefore,
\[
\frac 1 {x^2}  \left ( \int_{0}^{Q(v)}    \left (  \int_v^{H(t)} \left( \alpha^{-1}_{1,{\bf
X}} (u) \wedge n \right ) du  \right )^{1/2} dt   \right )^2 \leq  \frac{2}{x^{\eta }(2-\eta)}  \int_v^1  
\left( R_n^{-1}(u) \right )^{ \eta -1}  Q(u) du \, , 
\]
which combined with \eqref{I2n1} gives
\begin{equation}   \label{I2n}
I_2(n)  \leq   \frac{16 n}{x} \int_0^v Q(u) du  + \frac{32n}{x^{\eta }(2-\eta)}   \int_v^1 
\left ( R_n^{-1}(u) \right )^{ \eta -1}  Q(u) du 
\, .
\end{equation}
With similar arguments, we get
\begin{equation}   \label{I3n}
I_3(n)  \leq  \frac{16n}{x}  \int_0^v Q(u) du  + \frac{32n}{x^{\eta }(2-\eta)} \int_v^1  
\left( R_n^{-1}(u) \right)^{ \eta -1}  Q(u) du 
\, .
\end{equation}
Starting from  \eqref{decavecmartproba} and using  the upper bounds \eqref{I1n}, \eqref{I2n} and \eqref{I3n}, we derive that
\begin{multline*} 
{\mathbb P} \left ( \sup_{f \in \Lambda_1} \sum_{i=1}^n \left (f(X'_i) - \E(f(X'_i))\right ) \geq 5 x \right )   \\
\leq
\frac{34 n}{x} \int_0^v Q(u) du  + \frac{64 n}{x^{\eta } (2-\eta) } \int_v^1  ( R_n^{-1}(u) )^{ \eta -1}  Q(u) du   \, ,
\end{multline*}
which combined with \eqref{dec14FN} ends the proof of the proposition. $\diamond$

\subsection{A Rosenthal type inequality}\label{Sec:Ros}
In this section, we give some upper bounds for the 
quantity $\|W_1(\mu_n, \mu)\|_p$ when $p >2$ 
 in terms of the 
coefficients $\alpha_{2, \bf X}(k)$ and of the quantile function $Q$.  The function $\alpha_{2, {\bf X}}^{-1}$
is defined as in \eref{alphamoins1} by replacing  the coefficient 
$\alpha_{1, \bf X}(k)$ by
$\alpha_{2, \bf X}(k)$. 

\begin{prop}\label{propRos}
For $p>2$,  the following inequality holds:
\beq \label{RosenthalineW1}
\Vert  W_1 (\mu_n , \mu) \Vert_p^p \ll  \frac{s_{\alpha,n}^{p}}{n^{p/2}} + \frac 1 {n^{p-1}} \int_0^1 \left (\alpha_{2, {\bf X}}^{-1} (u) \wedge n \right )^{p-1}   Q^{p}(u) du \, ,
\eeq
where
\[
s_{\alpha,n}= \int_0^{\infty} \sqrt{S_{\alpha,n} (t)} dt 
\]
and  $S_{\alpha,n}$ is the function defined in 
\eref{def:Sn}.
\end{prop}

Note that Inequality \eref{RosenthalineW1} writes also
$$
\|  W_1 (\mu_n , \mu) \Vert_p^p \ll 
\frac{s_{\alpha,n}^{p}}{n^{p/2}} + \frac 1 {n^{p-1}}
  \sum_{k=0}^n (k+1)^{p-2}\int_0^{\alpha_{1, \bf X}(k)} Q^{p}(u) du \, .
$$

\begin{rem}
Inequality  \eref{RosenthalineW1} is similar to the Rosenthal  inequality 
for partial sums given in Theorem 6.3 
of Rio \cite{Ri}, with however two main differences:
\begin{itemize}
\item Firstly, the variance terms is not
the same, but this is because we consider the quantity $W_1 (\mu_n , \mu)$
and not only the partial sums, in accordance with the upper bounds for 
$\|  W_1 (\mu_n , \mu) \Vert_2$ given in Subsection \ref{M1/2}. 
\item Secondly, Rio's inequality is stated for 
$\alpha$-mixing sequences in the sense of Rosenblatt \cite{R}, and its proof 
relies on the coupling properties of these coefficients. 
Our result is valid for the larger class of $\alpha$-dependent sequences 
as defined in  \ref{defalpha} (with $k=2$ for the index of the dependency), 
and the proof is based on a version of the Rosenthal inequality for martingales given 
in \cite{MP}. Note that Rio's inequality cannot be applied 
to GPM maps, because the associated Markov chain is not $\alpha$-mixing  in the sense of 
Rosenblatt.
\end{itemize}

\end{rem}

\begin{rem} Let $r \geq 1$ and $p>2$. Starting again from \eref{Ebr}
and following the proof of Proposition \ref{propRos}, we obtain the 
upper bound
\beq \label{RosenthalineWr}
\Vert  W_r^r (\mu_n , \mu) \Vert_p^p \ll \frac{1}{n^{p/2}} \left (\int_0^{\infty} t^{r-1}
\sqrt{S_{\alpha,n} (t)} dt  \right)^p + \frac 1 {n^{p-1}} \int_0^1
 \left (\alpha_{2, {\bf X}}^{-1} (u) \wedge n \right)^{p-1}   Q^{rp}(u) du \, .
\eeq
\end{rem}

\begin{rem} Inequality \eqref{RosenthalineW1} implies in particular that if $p >2$ and 
\beq \label{condp>2}
\int_0^1 \left ( \alpha_{2,{\bf X}}^{-1}(u)\right)^{p/2}  Q^p (u) du < \infty \, ,
\eeq
then
$$
\Vert W_1 (\mu_n , \mu) \Vert_p \ll \frac 1 {\sqrt n} \, .
$$
\end{rem}

\begin{rem}
In the $m$-dependent case, Inequality \eref{RosenthalineW1}
becomes
$$
\|  W_1 (\mu_n , \mu) \|_p^p \ll 
\frac{1}{n^{p/2}} \left(\int_0^\infty \sqrt{H(t)} dt \right)^p 
+ \frac{1}{n^{p-1}} \|X_0\|_p^p \, .
$$
This inequality seems to be new even in the i.i.d. case. 
Compared to the usual Rosenthal bound for sums of i.i.d. random 
variables, the variance term is replaced by the integral 
involving $H$, in accordance  with the upper bound
\eref{M2}.
\end{rem}

\noindent {\bf Example (continued).}  We continue the example of Subsection \ref{CLT}. 
\begin{enumerate}
\item Let $p >2$, and let  $g$ be positive and non increasing on (0, 1), with
\[
g(x) \leq \frac{C}{x^{b}}
\quad \text{near 0, for some $C>0$ and $b \in [0, (1-\gamma)/p)$.}
\]
 Applying Proposition  \ref{propRos}, the following upper bounds hold.

For  $\gamma \in (0, 1/2)$ 
\beq 
\|W_1(\tilde \mu_n, \mu))\|_p
\ll
\begin{cases}
   n^{-1/2} \quad \quad \quad  \ \ \, \,  \text{if $b\leq (2-\gamma(p+2))/2p$}\\
   n^{(pb+\gamma-1)/p \gamma}  \quad  
  \text{if $ b >(2-\gamma(p+2))/2p $.}
  \end{cases}
  \nonumber
\eeq

For $\gamma \in [1/2, 1)$, 
$
\|W_1(\tilde \mu_n, \mu))\|_p 
\ll  n^{(pb+\gamma-1)/p \gamma}
$.
\item Let $p \in (0,1)$, and let  $g$ be positive and non decreasing on (0, 1), with
\[
g(x) \leq \frac{C}{(1-x)^{b}}
\quad \text{near 1, for some $C>0$ and $b \in [0, 1/p)$.}
\]
 
  Applying Proposition  \ref{propRos}, the following upper bounds hold.

For  $\gamma \in (0, 1/2)$ 
\beq 
\|W_1(\tilde \mu_n, \mu))\|_p
\ll
\begin{cases}
   n^{-1/2} \quad \quad \quad  \quad \ \ \, \,  \text{if $b\leq (2-\gamma(p+2))/2p(1-\gamma)$}\\
    n^{(\gamma-1)(1-pb)/p \gamma}  \quad  
  \text{if $ b >(2-\gamma(p+2))/2p(1-\gamma)$.}
  \end{cases}
  \nonumber
\eeq

For $\gamma \in [1/2, 1)$, 
$
\|W_1(\tilde \mu_n, \mu))\|_p 
\ll  n^{(\gamma-1)(1-pb)/p \gamma}
$.
\end{enumerate}

\begin{rem} In the case where $\theta$ is the LSV map defined 
by \eref{LSVmap} and $g$ is the identity (which is a particular
case of  Item 2, $b=0$, of the example above) all the rates 
for $\|W_1(\tilde \mu_n, \mu))\|_p $ given in Subsections  \ref{M1/2}, \ref{Sec:VBE} and 
\ref{Sec:Ros}
have been obtained in Corollary 4.1 of \cite{DM}  by  using a different approach. Moreover,   all the bounds are   
 optimal in that case
(see the discussion in Section 4.2 of \cite{DM}). 
\end{rem}

\medskip

\noindent {\bf Proof of Proposition \ref{propRos}}. Inequality \eqref{RosenthalineW1} follows from  Proposition \ref{propvineRosprob} below. 
\begin{prop}\label{propvineRosprob}
There exists a positive universal constant $c$ such that, for any positive integer $n$, any  $x >0$, any  $\eta >2$ and any  $\beta \in (\eta -2, \eta)$,  the following inequality holds:
\begin{multline} \label{ineRosprob}
{\mathbb P} \left ( n W_1 (\mu_n , \mu)  \geq   x \right ) \leq c  \frac{n^{\eta/2} }{x^{\eta}}  s_{\alpha,n}^{\eta} 
+  \frac{n}{x^{1+\beta/2}}  \int_0^{R_n^{-1}(x)}R^{\beta/2}_n(u)Q(u) du  \\+ c \frac{n}{x^{1+\eta/2}}  \int_{R_n^{-1}(x)}^1 R^{\eta/2}_n(u)Q(u) du   \, ,
\end{multline}
where \begin{equation*}
  R_n(u)=\left (\min\{q\in {\mathbb N}^* \tq \alpha_{2,{\bf X}}(q)\leq u\} \wedge n \right) Q(u)
 \  \text{ and } \ 
R_n^{-1}(x)=\inf \left \{u\in [0,1] \tq R_n(u)\leq x \right \}\,.
  \end{equation*}
\end{prop}
Indeed, 
\begin{multline} \label{p1rosmom}
\Vert n W_1 (\mu_n , \mu) \Vert_p^p =  p \int_0^{\infty} x^{p-1} {\mathbb P} \left ( n W_1 (\mu_n , \mu)  \geq x \right ) dx \\
\ll n^{p/2} s_{\alpha,n}^{p}+   \int_{n^{1/2}s_{\alpha,n}}^{\infty} x^{p-1} {\mathbb P} \left ( n W_1 (\mu_n , \mu)  \geq x \right ) dx \, .
\end{multline}
To handle the second term on the right-hand side, we apply \eqref{ineRosprob} with $\eta \in (2p-2, 2p)$ and $\beta \in (\eta -2  , 2p-2 )$. This gives 
\begin{multline*} 
  \int_{n^{1/2}s_{\alpha,n}}^{\infty} x^{p-1} {\mathbb P} \left ( n W_1 (\mu_n , \mu)  \geq x \right ) dx \ll n^{\eta/2}  s_{\alpha,n}^{\eta}  \int_{n^{1/2}s_{\alpha,n}}^{\infty} x^{p-\eta -1} dx 
  \\ + n  \int_0^{1}R^{\beta/2}_n(u)Q(u)   \int_{0}^{\infty}  x^{p-\beta/2-2}   {\bf 1}_{ u < R_n^{-1}(x)} du 
  \\ + n  \int_0^{1}R^{\eta/2}_n(u)Q(u)   \int_{0}^{\infty}  x^{p-\eta/2-2}   {\bf 1}_{ u \geq R_n^{-1}(x)} du \, .
\end{multline*} 
Since  $ u < R_n^{-1}(x) \iff x < R_n(u)$, the choice of $\eta$ and $\beta$ implies that, for any $p>2$,
\[
   \int_{n^{1/2}s_{\alpha,n}}^{\infty} x^{p-1} {\mathbb P} \left ( n W_1 (\mu_n , \mu)  \geq x \right ) dx \ll n^{p/2} s_{\alpha,n}^{p} + n \int_0^1 R_n^{p-1}(u)   Q(u) du  \, ,
\]
which together with \eqref{p1rosmom} give \eqref{RosenthalineW1}. 

 To complete  the proof of Proposition \ref{propRos}, it remains to prove Proposition \ref{propvineRosprob}. With this aim, we proceed as for the proof of Proposition \ref{propvBEine} with the following modification: in the definition of $R_n$ (and then also of $v$ defined in \eqref{defvM}), $\alpha_{1,{\bf
X}}$ is replaced by $\alpha_{2,{\bf
X}}$, and in the definition of $q$ given in \eqref{defofq}, $\alpha_{1,{\bf
X}}$ is also replaced by $\alpha_{2,{\bf
X}}$. Assuming first that $q=n$, we first notice,  by following the proof of Proposition \ref{propvBEine}, that the bound \eqref{dec13FNdef} is still valid. In addition since $ u < R_n^{-1}(x) \iff x < R_n(u)$,
\beq \label{p1Rosprob}
\int_0^{R^{-1}_n(x)} Q(u) du \leq x^{-\beta/2} \int_0^{R^{-1}_n(x)} R^{\beta/2}_n(u)Q(u) du \, ,
\eeq
which combined with \eqref{dec13FNdef}   proves the proposition in case where $q=n$.

From now on, we assume that $q<n$ (therefore $\alpha_{2,{\bf
X}} (q)\leq v$). The bound \eqref{dec14FN} is still valid and combined with \eqref{p1Rosprob} gives 
\begin{multline}\label{p2Rosprob}
   {\mathbb P} \left ( n W_1 (\mu_n , \mu)  \geq 6x \right )  
  \leq   {\mathbb P} \left ( \sup_{f \in \Lambda_1} \sum_{i=1}^n \left ( f(X'_i) - \E(f(X'_i))\right )
   \geq 5 x \right)  \\ +  \frac{2n}{x^{1+\beta/2}}  \int_0^{R^{-1}_n(x)} R^{\beta/2}_n(u)Q(u) du \, .
 \end{multline}
As in the proof of Proposition \ref{propvBEine}, the first term on the right-hand side can be handled with the help of the decomposition \eqref{decavecmartproba}.  Clearly since $\alpha_{1,{\bf
X}} (q) \leq \alpha_{2,{\bf
X}} (q) \leq v$, the term $I_1(n)$ in \eqref{decavecmartproba} satisfies the inequality \eqref{I1n}. Therefore taking into account \eqref{p1Rosprob}, it follows that 
\begin{equation} \label{p3I1n}
I_1(n) 
\leq \frac{2n}{x^{1+\beta/2}}  \int_0^{R^{-1}_n(x)} R^{\beta/2}_n(u)Q(u) du   \, .
\end{equation}
We handle now the term $I_2(n)$ in the decomposition \eqref{decavecmartproba}. Using again Markov's inequality, we get that for any $\eta >2$,
\[
I_2(n)  \leq \frac 1 { x^{\eta} }\left ( \int_{-Q(v)}^{Q(v)}   \left \Vert \max_{2 \leq 2j \leq
 [n/q]
 } \left \vert \sum_{i=1}^j \tilde U_{2i} (t) \right \vert \right \Vert_{\eta} dt   \right )^{\eta} \, .
\]
Note that $( \tilde U_{2i}(t))_{i \in {\mathbb Z}} $ (resp. $( \tilde U_{2i -1
}(t))_{i \in {\mathbb Z}}$) is a stationary sequence of martingale differences  with respect to the filtration $ ( {\cal {G}}_{2i} )_{i \in {\mathbb Z}}$ (resp. $ ( {\cal {G}}_{2i -1 } )_{i \in {\mathbb Z}}$). By using the Rosenthal inequality of Merlev\`ede and Peligrad \cite{MP} for martingales (see their Theorem 6), we get
\begin{multline*}
\left \Vert \max_{2 \leq 2j \leq
 [n/q]
 } \left \vert \sum_{i=1}^j \tilde U_{2i} (t) \right\vert \right\Vert_\eta \\
\ll (n/q)^{1/\eta}  \left \Vert \tilde U_{2} (t)  \right \Vert_\eta + (n/q)^{1/\eta} \left (  \sum_{k=1}^{[n/q]} \frac{1}{k^{1+2\delta/\eta}}  \left \Vert  \E_0 \left ( \left ( \sum_{i=1}^k \tilde U_{2i} (t) \right )^2
\right )\right \Vert_{\eta/2}^{\delta}\right )^{1/(2\delta)}\, ,
\end{multline*}
where $\delta= \min \left \{1, (\eta-2)^{-1}\right \}$. Since $( \tilde U_{2i}(t))_{i \in {\mathbb Z}} $ is a stationary sequence of martingale differences  with respect to the filtration $ ( {\cal {G}}_{2i} )_{i \in {\mathbb Z}}$,
\[
\E_0 \left ( \left ( \sum_{i=1}^k \tilde U_{2i} (t) \right )^2 \right ) = \sum_{i=1}^k \E_0 
\left (  \tilde U^2_{2i} (t)   \right ) \, .
\]
Moreover $
\E_0 \big (  \tilde U^2_{2i} (t)   \big ) \leq \E_0 \left ( U^2_{2i} (t)   \right ) $. 
Therefore
\[
\left \Vert  \E_0 \left ( \left ( \sum_{i=1}^k \tilde U_{2i} (t) \right )^2\right )\right \Vert_{r/2} \leq 
\sum_{i=1}^k \left \Vert \E_0 \left ( U^2_{2i} (t) \right ) - \E \left (   U^2_{2i} (t)  \right )\right \Vert_{r/2} +  \sum_{i=1}^k  \E \left (U^2_{2i} (t) \right ) \, .
\]
By stationarity
\[
\sum_{i=1}^k  \E \left ( U^2_{2i} (t) \right ) = k \left \Vert S'_q(t)  \right \Vert_2^2 \, ,
\]
where
\[
S_q'(t) = \sum_{i=1}^q \left ( \I_{X_i' \leq t} - \E \left ( \I_{X_i' \leq t} \right ) \right )  \, .
\]
It follows that
\begin{multline*}
\left\Vert \max_{2 \leq 2j \leq
 [n/q]
 } \left \vert \sum_{i=1}^j \tilde U_{2i} (t) \right\vert \right \Vert_r \\
\ll (n/q)^{1/\eta}  \left \Vert  S'_q (t) \right \Vert_\eta + (n/q)^{1/2}  \left \Vert S'_q (t) \right \Vert_2 +  (n/q)^{1/\eta} \left  (  \sum_{k=1}^{[n/q]} \frac{1}{k^{1+2\delta/r}} D_{k,q}^{\delta}(t)
\right )^{1/(2\delta)}\, ,
\end{multline*}
where
\[
D_{k,q} (t) = \sum_{i=1}^k  \left \Vert \E_0 \left ( U^2_{2i} (t) \right ) - \E \left (   U^2_{2i} (t)  \right )\right \Vert_{\eta/2}\, .
\]
We have
\begin{align*}
D_{k,q} (t) 
 \leq & q^2\sum_{i=1}^k   \sup_{j \geq \ell \geq (i-1)q+1}\sup_{t \in {\mathbb R}}  \big \Vert \E_0 \big ((\I_{X_\ell' \leq t} - {\mathbb E} (\I_{X_\ell' \leq t} ) ) (\I_{X_j \leq t}  - {\mathbb E} (\I_{X_j' \leq t} ) )    \big ) \\
&  \hspace{5cm}- \E\big ((\I_{X_\ell' \leq t} - {\mathbb E} (\I_{X_\ell' \leq t} ) ) (\I_{X_j' \leq t} - {\mathbb E} (\I_{X_j' \leq t} ) )    \big ) \big \Vert_{\eta/2} \\
\leq & q^2 \sum_{i=1}^k   \sup_{j \geq \ell \geq (i-1)q+1}\sup_{t \in {\mathbb R}}  \big \Vert \E_0 \big ((\I_{X_\ell' \leq t} - {\mathbb E} (\I_{X_\ell' \leq t} ) ) (\I_{X_j \leq t}  - {\mathbb E} (\I_{X_j' \leq t} ) )    \big ) \\
& \hspace{5cm} - \E\big ((\I_{X_\ell' \leq t} - {\mathbb E} (\I_{X_\ell' \leq t} ) ) (\I_{X_j' \leq t} - {\mathbb E} (\I_{X_j' \leq t} ) )    \big ) \big \Vert_{\eta/2} \\
\leq  & q^2 \sum_{i=1}^k  \alpha_{2,{\bf
X}}^{2/\eta} (iq+1) 
\, ,
\end{align*}
where we have used the fact that $g_M$ is nondecreasing for the second inequality. Since $\beta < \eta$,  H\"older's inequality gives
\[
D_{k,q} (t) \ll q^2 k^{ (\eta  - \beta)/\eta} \left ( \sum_{i=1}^k i^{\beta/2-1}\alpha_{2,{\bf
X}} (iq+1)  \right )^{2/\eta} \, .
\]
Therefore, since $\beta > \eta -2$,
\begin{multline*} \left ( \int_{-Q(v)}^{Q(v)}   \frac{n^{1/\eta}}{q^{1/\eta}} \left (  \sum_{k=1}^{[n/q]} \frac{1}{k^{1+2\delta/r}} D_{k,q} (t)^{\delta}\right )^{1/(2\delta)}  dt   \right )^\eta  \\
\ll  n q^{\eta - 1 }    Q^\eta(v) \left (  \sum_{k=1}^{[n/q]} \frac{k^{\delta(\eta  - \beta)/\eta} }{k^{1+2\delta/\eta}} \right )^{\eta/(2\delta)}    \sum_{i=1}^{[n/q]} i^{\beta/2-1} \alpha_{2,{\bf
X}}(iq+1) \\
\ll  n q^{\eta - 1 }   Q^\eta(v)     \sum_{i=1}^{[n/q]} i^{\beta/2-1} \alpha_{2,{\bf
X}} (iq+1)  \, .
 \end{multline*}
Note that since $y < \alpha_{2,{\bf
X}} ^{-1} (u) \iff \alpha_{2,{\bf
X}} (y) >u$ and $\alpha_{2,{\bf
X}} (q)\leq v$,
 \begin{multline*}
 \sum_{i=1}^{[n/q]} i^{\beta/2-1} \alpha_{2,{\bf
X}} (iq+1)  =  \sum_{i=1}^{[n/q]} i^{\beta/2-1} \int_0^1 \I_{u < \alpha _{2,{\bf
X}} (iq+1)} \\
 \leq \int_0^v \sum_{i=1}^{[n/q]} i^{\beta/2-1} \I_{ i \leq q^{-1} \alpha_{2,{\bf
X}} ^{-1} (u)}
 \leq   q^{-\beta/2}\int_0^v \left ( \alpha_{2,{\bf
X}} ^{-1} (u) \wedge n\right)^{\beta/2} du \, .
\end{multline*}
Hence
\begin{equation*}\left ( \int_{-Q(v)}^{Q(v)}   \frac{n^{1/\eta}}{q^{1/\eta}} \left (  \sum_{k=1}^{[n/q]} \frac{1}{k^{1+2\delta/r}} D_{k,q} (t)^{\delta}\right )^{1/(2\delta)}  dt   \right )^\eta 
\ll n   q^{\eta -1-\beta/2} Q^\eta(v) \int_0^v \left ( \alpha_{2,{\bf
X}} ^{-1} (u) \wedge n \right )^{\beta/2} du  \, .
 \end{equation*}
Using \eqref{restq1} and the fact that $u<v \iff Q(v) <Q(u)$, we infer that
\begin{equation} \label{majDkdelta}
\left ( \int_{-Q(v)}^{Q(v)}   \frac{n^{1/\eta}}{q^{1/\eta}} \left (  \sum_{k=1}^{[n/q]} \frac{1}{k^{1+2\delta/r}} D_{k,q} (t)^{\delta}\right )^{1/(2\delta)}  dt   \right )^\eta 
\ll n x^{\eta-\beta/2-1}   \int_0^v R_n^{\beta/2} (u)  Q(u) du  \, .
 \end{equation}
 On another hand, since 
  \begin{equation*}
 \left \Vert S'_q (t) \right \Vert_2^2 \leq 2q  \sum_{k=0}^{q-1} \left | \E \left ( (\I_{g_M(X_0) \leq t}- {\mathbb E} ( \I_{g_M(X_0) \leq  t} ) )  (\I_{g_M(X_k) \leq  t}- {\mathbb E} ( \I_{g_M(X_k)  \leq  t} ) ) \right )  \right |\, ,
   \end{equation*}
   proceeding as to bound $I_2(n)$ in the proof of Proposition \ref{propvBEine}, we infer that 
 \beq\label{butSq2}
  \left ( \int_{-Q(v)}^{Q(v)}   \frac{n^{1/2}}{q^{1/2}} \left \Vert S'_q (t)  \right \Vert_2 dt   \right)^\eta  \ll n^{\eta/2} s_{\alpha,n}^\eta  \, .
 \eeq
 We prove now that
 \beq\label{butSqr}
 \left ( \int_{-Q(v)}^{Q(v)}   \frac{n^{1/\eta}}{q^{1/\eta}} \left \Vert S'_q (t) \right 
 \Vert_\eta dt   \right )^\eta  \ll n x^{\eta-\beta/2-1}   \int_0^v R_n^{\beta/2} (u)  Q(u) du  
 + n x^{\eta/2-1}  \int_{v}^1 R^{\eta/2}_n(u)Q(u) du   \, .
 \eeq
With this aim, assume first that we can prove that
\beq \label{MarcSq}
\int_{-Q(v)}^{Q(v)}  \left \Vert S'_q (t) \right \Vert_\eta dt \ll q^{1/2} \int_{0}^{Q(v)} \left ( \int_0^{H(t)} \left (\alpha_{2,{\bf
X}}^{-1}(u) \wedge q\right )^{\eta /2}du\right )^{1/\eta} dt \, ,
\eeq
then 
\[
 \left ( \int_{-Q(v)}^{Q(v)}   \frac{n^{1/\eta}}{q^{1/\eta}} \left \Vert S'_q (t)  \right \Vert_\eta dt   \right
  )^\eta \ll A(n) + B(n) \, ,
\]
where
\[
A(n) =  n v q^{\eta-1}  Q^\eta(v) 
\quad 
\text{and}
\quad 
B(n) =  n  q^{\eta/2-1}  \left (\int_{0}^{Q(v)} \left ( \int_v^{H(t)} \left (\alpha_{2,{\bf
X}}^{-1}(u)  \wedge q \right )^{\eta/2}du\right )^{1/\eta} dt   \right )^\eta  dx \, .
\]
Using \eqref{restq1}, the fact that $u<v \iff Q(v) <Q(u)$ and that $u < R_n^{-1} (x)=v \iff x < R_n(u)$, we successively derive
\beq \label{ubAn}
A(n) \ll n  x^{\eta-1}  v Q(v) \ll  n  x^{\eta-1}    \int_0^v Q(u) dx  \ll n x^{\eta-\beta/2-1}   \int_0^v R_n^{\beta/2} (u)  Q(u) du \, .
\eeq
On the other hand,  since  $u < H(t) \iff t < Q(u)$, we have
\begin{multline*}
 B(n) \leq  n q^{\eta/2-1}  \left (\int_{0}^{Q(v)} \frac{1}{t^{1/2+1/\eta}}\left ( \int_v^{H(t)} 
 \left ( \alpha_{2,{\bf
X}}^{-1}(u) \wedge q\right )^{r/2} Q^{\eta/2+1}(u)du\right  )^{1/\eta} dt   \right  )^\eta dx \\
 \ll n( q Q(v) )^{\eta/2-1}   \int_v^{1} \left (\alpha_{2,{\bf
X}}^{-1}(u) \wedge n \right )^{\eta /2} Q^{\eta/2+1}(u)du \, .
  \end{multline*}
Using \eqref{restq1}, it follows that
\[ B(n)
 \ll n x^{\eta/2-1}   \int_v^{1} R_n^{\eta /2} (u) Q(u)du   \, .
\]
This last upper bound together with \eqref{ubAn} show that to prove \eqref{butSqr} it suffices to prove \eqref{MarcSq}. To prove this moment inequality, we use Corollary 2 in \cite{DD}. Since, for any $t \in {\mathbb R}$, $\vert \I_{X_0' \leq t} - {\mathbb E} (\I_{X_0' \leq t} )  \vert \leq 1$, this gives
\[
 \left \Vert S'_q (t)  \right \Vert_\eta \leq \sqrt{2 q \eta } \left (  \int_0^{\Vert Y(t) \Vert_1} 
 \left (\gamma^{-1} (u) \wedge q \right )^{\eta/2} du \right )^{1/\eta} \, ,
\]
where $Y(t) = \I_{g_M(X_0) \leq t} - {\mathbb E} (\I_{g_M(X_0) \leq t} ) $ and 
$$
\gamma^{-1} (u)  = \sum_{k=0}^{\infty}  \I_{u \leq \gamma(k)}  \quad  \text{with } \quad 
 \gamma(k) = \left \Vert \E_0 \left (  \I_{g_M(X_k) \leq  t} - {\mathbb E} \left  (\I_{g_M(X_k) \leq t} 
 \right )  \right )   \right \Vert _1 \, .
$$
Since $g_M$ is nondecreasing 
$
\gamma(k) 
\leq  \alpha_{1,{\bf
X}}(k) \leq \alpha_{2,{\bf
X}}(k)  $ in such a way that
 $\gamma^{-1} (u)  \leq  \alpha_{2,{\bf
X}}^{-1}(u) $. 
Moreover, for any $t \in {\mathbb R}$,  
$$
\Vert Y(t) \Vert_1 = 2 {\mathbb P}(g_M(X_0) \leq  t){\mathbb P}(g_M(X_0) >t   )
\leq 2 \min \left \{ {\mathbb P}(|X_0| \geq  - t), {\mathbb P}(|X_0| >  t) \right \} \, .
$$
All these considerations end the proof of  \eqref{MarcSq}. 

So, overall, we get 
\begin{equation*} \label{p3I2n}
I_2(n) 
\ll x^{- \eta }n^{\eta/2} s_{\alpha,n}^\eta  + n x^{-\beta/2-1}   \int_0^v R_n^{\beta/2} (u)  Q(u) du  + n x^{-\eta/2-1}  \int_{v}^1 R^{\eta/2}_n(u)Q(u) du    \, .
\end{equation*}
With similar arguments, we can prove that 
\begin{equation*} \label{p3I3n}
I_3(n) 
\ll x^{- \eta }n^{\eta/2} s_{\alpha,n}^\eta  + n x^{-\beta/2-1}   \int_0^v R_n^{\beta/2} (u)  Q(u) du  + n x^{-\eta/2-1}  \int_{v}^1 R^{\eta/2}_n(u)Q(u) du    \, .
\end{equation*}
Therefore starting from \eqref{decavecmartproba}  and taking into account \eqref{p3I1n}, \eqref{p3I2n} and \eqref{p3I3n}, it follows that 
\begin{multline*} 
{\mathbb P} \left ( \sup_{f \in \Lambda_1} \sum_{i=1}^n \left ( f(X'_i) - \E(f(X'_i))\right ) \geq 5 x \right )   \\
\ll
x^{- \eta }n^{\eta/2} s_{\alpha,n}^\eta  + n x^{-\beta/2-1}   \int_0^v R_n^{\beta/2} (u)  Q(u) du  + n x^{-\eta/2-1}  \int_{v}^1 R^{\eta/2}_n(u)Q(u) du     \, ,
\end{multline*}
which combined with \eref{p2Rosprob} ends the proof of Proposition \ref{propvineRosprob}. $\diamond$

\section{Weak convergence of partial sums in ${\mathbb L}^1(m)$}\label{CLTL1}

\setcounter{equation}{0}

Let $(S, {\mathcal S}, m)$ be a $\sigma$-finite measure 
space such that 
${\mathbb L}^1(S, {\mathcal S}, m)$ is separable. In what follows, we shall 
denote by 
${\mathbb L}^1(m)$ the space 
${\mathbb L}^1(S, {\mathcal S}, m)$. 

We use the notations of Section 
\ref{DefNot}. Let $Y_0=\{Y_0(t), t \in S\}$ be a random variable with values 
in ${\mathbb L}^1(m)$, such that
$$
\int \|Y_0(t)\|_1  \ m(dt) < \infty \quad \text{and} \quad 
\int Y_0(t)   \ m(dt) =0 \, .
$$
Define the stationary sequence 
 ${\bf Y}=(Y_i)_{i \in {\mathbb Z}}$ by $Y_i=Y_0 \circ T^i$, and let 
 $$
 S_n=\sum_{k=1}^n Y_k  \, . 
 $$

\subsection{Previous results}\label{sec:CD}

If ${\bf Y}$ is  a sequence of i.i.d. random variables,  Jain \cite{J} proved that 
$n^{-1/2} S_n$ 
satisfies the CLT (i.e. converges in distribution to an
${\mathbb L}^1(m)$-valued Gaussian random variable) if and only if 
\begin{equation}\label{Jain}
\int \|Y_0(t)\|_2 \ m(dt) < \infty. 
\end{equation}

Using a general result by de Acosta, Araujo and Gin\'e
\cite{AAG}, D\'ed\'e \cite{D} proved that the CLT 
remains valid under \eref{Jain} for stationary and ergodic martingale differences (meaning that 
${\mathbb E}(Y_1|{\mathcal F}_0)=0$ almost surely).
Starting from a martingale approximation, she proved then that, if $\bf Y$ is ergodic, 
the CLT holds  as soon as
\eref{Jain}  holds and 
\begin{equation}\label{Dede}
\sum_{k \in {\mathbb Z}} \int \|P_0(Y_k (t))\|_2 \ m(dt) < \infty , 
\end{equation}
where 
${\mathbb P}_0(Y_k(t))={\mathbb E}(Y_k(t)|{\mathcal F}_0)
-{\mathbb E}(Y_k(t)|{\mathcal F}_{-1})$. 

In a recent paper, Cuny \cite{C} has given many new results concerning the behavior of partial sums of
dependent sequences in Banach spaces of cotype 2. 
Among these results, he showed that, if ${\bf Y}$ is ergodic, 
$Y_0$ is ${\mathcal F}_0$-measurable,  \eref{Jain}  holds
and 
\begin{equation}\label{Cuny}
\sum_{n>0} \int 
\frac{\|{\mathbb E}(S_n|{\mathcal F}_0)\|_2}{n^{3/2}} \ m(dt) < 
\infty \, , 
\end{equation}
then the CLT and the weak invariance principle (WIP) hold.
By WIP, we mean that the partial sum process
$
\{n^{-1/2} S_{[nt]}, t \in [0,1]\}
$
converges in distribution to an 
${\mathbb L}^1(m)$-valued Wiener process in the space $D_{{\mathbb L}_1(m)}([0,1])$
of ${\mathbb L}^1(m)$-valued c\`adl\`ag functions  equipped with the uniform metric. As usual, an ${\mathbb L}^1(m)$-valued Wiener process  with covariance 
$\Lambda$ is   a  centered Gaussian process 
$W=\{W(t), t \in [0,1]\}$ such that
$\mathbb E(\|W(t)\|_{{\mathbb L}^1(m)}^2< \infty$ 
for all $t\in [0,1]$
and, for all $f,g$ in ${\mathbb L}^\infty(m)$, 
$$
{\mathrm {Cov}}  \left ( \int f(u) W_t(u) \ m(du), 
\int g(u) W_s(u) \  m(du) \right) = \min \{s,t\} \Lambda(f,g)
$$
(as usual, we identify a function $f$ in ${\mathbb L}^\infty(m)$ with an
element of the dual of ${\mathbb L}^1(m)$). 

Note that Cuny \cite{C} also proved that the WIP
holds under \eref{Dede}, and that the almost sure invariance principle with rate $o(\sqrt {n \ln \ln n})$ is 
true if either \eref{Dede} of \eref{Cuny} holds. 

The condition \eref{Dede} is the ${\mathbb L}^1(m)$ version of Hannan's criterion \cite{H}, and the condition 
\eref{Cuny}   is the ${\mathbb L}^1(m)$ version of 
Maxwell-Woodroofe's criterion criterion \cite{MW}. 
If $Y_0$ is ${\mathcal F}_0$-measurable, 
both criteria hold as soon as 
\begin{equation}\label{conseq}
\sum_{k=0}^\infty \frac{1}{k+1} \int  \|{\mathbb E}(
Y_k(t)|{\mathcal F}_0)\|_2  \ m(dt) < \infty \, . 
\end{equation}

As shown  in \cite{C}, if either 
\eref{Dede} or \eref{Cuny} holds, there exists
a stationary and  ergodic sequence of martingale 
differences $(D_i)_{i \in {\mathbb Z}}$ with values in
${\mathbb L}^1(m)$, such that, setting $M_n=\sum_{k=1}^n D_k$,  
$$
  \left \| \max_{1 \leq k \leq n} 
  \int \left | S_k(t)- M_k(t)\right | \ m(dt) \right \|_2 = o( \sqrt n) \, . 
$$
In the next subsections, we shall rather look 
for a martingale approximation in ${\mathbb L}^1$, 
in the spirit of Gordin \cite{G}. Our criterion will not
be directly comparable to  either 
\eref{Dede} or \eref{Cuny}, but its application to the empirical distribution function of $\alpha$-dependent
sequences will lead to weaker conditions (see Section 
\ref{quantile}
for a deeper discussion). 

\subsection{A central limit theorem in ${\mathbb L}^1(m)$ for non-adapted sequences}

In this section,  we give an extension of  Gordin's criterion  \cite{G} for
the central limit theorem to ${\mathbb L}^1(m)$-valued
random variables.

\begin{thm}\label{thm:nonadapted}
Assume that, for $m$ almost every $t$, the series
\begin{equation}\label{GL1}
 U(t)=\sum_{k=1}^\infty {\mathbb E}_0(Y_k(t)) \quad  \text{and} \quad  V(t)=-\sum_{k=-\infty}^0 \big (Y_{k}(t) - {\mathbb E}_0(Y_{k}(t)) \big)
\end{equation}
converge in probability, and let
$$
D_0(t)= \sum_{k \in {\mathbb Z}} \big( {\mathbb E}_0(Y_k(t))- {\mathbb E}_{-1}(Y_k(t))\big)
\quad \text{and} \quad M_n(t)= \sum_{k=1}^n D_0(t) \circ T^k \, .
$$
If
\begin{equation}\label{UVcond}
\int \|U(t)+V(t)\|_1 \ m(dt) < \infty \,  ,
\end{equation}
then
\begin{equation}\label{approx}
\lim_{n \rightarrow \infty} \int \left \| \frac{S_n (t)}{\sqrt n}- \frac{M_n(t)}{\sqrt n} 
\right \|_1 \ m (dt)=0 \, .
\end{equation}
If moreover, for $m$ almost every $t$,
\begin{equation}\label{condGordin}
 C(t)=\liminf_{n \rightarrow \infty} \frac{1}{\sqrt n} {\mathbb E}(|S_n(t)|)< \infty \quad
 \text{and} \quad  \int C(t) \ m(dt) < \infty \, ,
\end{equation}
then
\begin{equation}\label{moment2}
 \int \|D_0(t)\|_2 \ m(dt) < \infty \, ,
\end{equation}
 and for any $(s_1, \ldots , s_d)$ in ${[0,1]}^d$, the random vector
$n^{-1/2} (S_{[ns_1]}, \ldots , S_{[ns_d]})^t$ converges in distribution
in $({\mathbb L}^1(m))^{d}$  to the Gaussian random vector $(W_{s_1}, \ldots , W_{s_d})$, where
$W$ is the ${\mathbb L}_1(m)$-valued Wiener process $W$
with covariance
operator $\Lambda$ defined by: for any $f,g$ in ${\mathbb L}_\infty(m)$,
\begin{equation}\label{defcov}
\Lambda(f,g)=  {\mathbb{E}}\left( \iint f(t)g(s) D_0(t) D_0(s) \ m(dt)m(ds) \right)\, .
\end{equation}
\end{thm}

\noindent {\bf Proof of Theorem \ref{thm:nonadapted}.} We first state the  following intermediate result:
\begin{prop}\label{prop:cobord}
Assume that, for $m$ almost every $t$,
\begin{equation}\label{cob}
  Y_0(t)= D_0(t) + Z(t) - Z(t) \circ T \, ,
\end{equation}
where $D_0(t)$ is an integrable random variable such that  ${\mathbb E}(D_0(t)|{\mathcal F}_{-1})=0$
almost surely.
Let then
$
M_n(t)= \sum_{k=1}^n D_0(t) \circ T^k \, .
$
If
\begin{equation}\label{Zcond}
\int \|Z(t)\|_1 \ m(dt) < \infty \, ,
\end{equation}
then \eref{approx} holds. If moreover \eref{moment2} holds, then the conclusion of
Theorem \ref{thm:nonadapted} holds.
\end{prop}
Before proving Proposition \ref{prop:cobord}, let us continue the proof of Theorem \ref{thm:nonadapted}.
Note first that, if \eref{GL1} is satisfied, then \eref{cob} holds, with
$$
D_0(t)= \sum_{k \in {\mathbb Z}} \big({\mathbb E}_0(Y_k(t))
-{\mathbb E}_{-1}(Y_k(t))\big) \quad \text{and} \quad
Z(t) \circ T= \sum_{k=1}^\infty {\mathbb E}_0(Y_k(t)) -\sum_{k=-\infty}^0 \big(Y_k(t)-{\mathbb E}_0(Y_k(t))\big)\, .
$$
Now, if $Z(t)$ is defined as above, the conditions \eref{UVcond} and \eref{Zcond} are the same.
Hence, it follows from Proposition \ref{prop:cobord} that \eref{approx} holds as soon as \eref{Zcond} is satified. The second part of Theorem \ref{thm:nonadapted} will follow from Proposition \ref{prop:cobord}
if we prove that \eref{condGordin} implies \eref{moment2}. By  \eref{UVcond} it follows that
\begin{equation}\label{liminfZ}
\lim_{n \rightarrow \infty} \frac{\|Z(t)\|_1}{\sqrt n}=0 \quad \text{for $m$-almost every $t$.}
\end{equation}
Since $S_n(t)=M_n(t) + Z(t) - Z(t)\circ T^n$, we infer from \eref{liminfZ} that,
for $m$ almost every $t$,
\begin{equation}\label{liminfM}
\liminf_{n \rightarrow \infty} \frac{\|M_n(t)\|_1}{\sqrt n}=
\liminf_{n \rightarrow \infty} \frac{\|S_n(t)\|_1}{\sqrt n} \, .
\end{equation}
From \eref{liminfM} and \eref{condGordin}, it follows that,
for $m$ almost every $t$,
$$
C(t)=\liminf_{n \rightarrow \infty} \frac{\|M_n(t)\|_1}{\sqrt n}< \infty \, .
$$
Now, applying Theorem 1 and Remark 1.1 in Esseen and Janson \cite{EJ}, we deduce that,
for
$m$ almost every $t$,
$$
\|D_0(t)\|_2= \sqrt{\frac \pi 2} C(t) \, ,
$$
so that \eref{condGordin} implies \eref{moment2}. This completes the proof of Theorem
\ref{thm:nonadapted}.

\medskip

\noindent {\bf Proof of Proposition \ref{prop:cobord}.}
Since $S_n(t)=M_n(t) + Z(t) - Z(t)\circ T^n$, it follows that
$$
\int \left \| \frac{S_n (t)}{\sqrt n}- \frac{M_n(t)}{\sqrt n} \right \|_1 \ 
m (dt) \leq \frac{2}{\sqrt n}
\int \|Z(t)\|_1 \ m(dt) \, ,
$$
and \eref{approx}
follows from \eref{Zcond}.

Now, let $d$ be a positive integer, and let $f$ be a separately Lipschitz function from
$({\mathbb L}^1(m))^d$ to ${\mathbb R}$. This means that there exists non-negative constants
$c_1, \ldots, c_d$ such that
$$
|f(x_1, \ldots, x_d)-f(y_1, \ldots , y_d)|\leq \sum_{i=1}^d c_i \int |x_i(t)-y_i(t)|  \ m(dt) \, .
$$
For such a $f$ and  any  $(s_1, \ldots, s_d)$ in $[0,1]^d$, we get that
\begin{multline*}
\left |{\mathbb E}\left ( f \left( \frac{S_{[ns_1]}}{\sqrt n}, \ldots , \frac{S_{[ns_d]}}{\sqrt n}\right) \right)-
{\mathbb E}\left ( f \left( \frac{M_{[ns_1]}}{\sqrt n}, \ldots , \frac{M_{[ns_d]}}{\sqrt n}
\right) \right) \right | \\
\leq \sum_{i=1}^d c_i \int \left \| \frac{S_{[ns_i]} (t)}{\sqrt n}- 
\frac{M_{[ns_i]}(t)}{\sqrt n} \right \|_1 \ m (dt)\, ,
\end{multline*}
and it follows from \eref{approx} that
\begin{equation}\label{reducM}
\lim_{n \rightarrow \infty} \left |{\mathbb E}\left ( f \left( \frac{S_{[ns_1]}}{\sqrt n}, \ldots , \frac{S_{[ns_d]}}{\sqrt n}\right) \right)-
{\mathbb E}\left ( f \left( \frac{M_{[ns_1]}}{\sqrt n}, \ldots , \frac{M_{[ns_d]}}{\sqrt n}\right)
 \right)\right |=0 \, .
\end{equation}
Now,  when \eref{moment2} holds, Cuny \cite{C}  proved that the process
$\{n^{-1/2} M_{[nt]}, t \in [0,1]\}$ converges in distribution
in  the space $D_{{\mathbb L}_1(m)}([0,1])$ to an ${\mathbb L}_1(m)$-valued Wiener process $W$, with covariance
operator $\Lambda$ given by \eref{defcov}. Together with \eref{reducM}, this completes the proof of Proposition \ref{prop:cobord}.

\subsection{An invariance principle in ${\mathbb L}^1(m)$ for adapted sequences}

In this subsection, we assume that the random variable
$Y_0$ is ${\mathcal F}_0$-measurable.
\begin{thm}\label{pifL1}
Assume that, for $m$-almost every $t$, the series 
$U(t)$ defined in \eref{GL1} converges in
probability. Assume also that, for $m$-almost every $t$, the series 
\begin{equation}\label{serieDR}
\sum_{k=0}^n Y_0(t) {\mathbb E}_0(Y_k(t)) 
\end{equation} converge in ${\mathbb L}^1$, and let 
\begin{equation}\label{Lt}
L(t)= \sup_{n \geq 0}
\left \|  \sum_{k=0}^n Y_0(t) {\mathbb E}_0(Y_k(t)) \right \|_1  \, .
\end{equation}
If moreover $\int \|U(t)\|_1 \ m(dt)  < \infty$  and
\begin{equation}
\label{condui}
 \quad \int \sqrt{L(t)} \ m(dt)  < \infty \, ,
\end{equation}
then $\{n^{-1/2} S_{[nt]}, t \in [0,1]\}$ converges in distribution
in  the space $D_{{\mathbb L}_1(m)}([0,1])$ to an ${\mathbb L}_1(m)$-valued Wiener process $W$, with covariance
operator $\Lambda$ defined by \eref{defcov}.
\end{thm}
As an immediate consequence of Theorem \ref{pifL1}, the following corollary holds:
\begin{cor}\label{easycor}
Assume that
\begin{equation}\label{easycrit}
\int \sqrt{\sum_{k\geq 0} \big \|  \max\{1,|Y_0(t)|\} |{\mathbb E}_0(Y_k(t))| \big \|_1}  \ m(dt) < \infty \, .
\end{equation}
Then the conclusion of Theorem \ref{pifL1} holds.
\end{cor}

\begin{rem}\label{rem:extent}
Under the assumptions of Theorem \ref{pifL1}, we shall prove that the sequence
\begin{equation}\label{maxui}
T_n= \frac 1 n \left ( \max_{1 \leq k \leq n} \int |S_k(t)| m(dt) \right)^2
\end{equation} 
is uniformly integrable (see Lemma \ref{ui} below). By standard arguments, this 
implies the following extension of Theorem \ref{pifL1}: 
let $\psi$ be any continuous  function 
from $(D_{{\mathbb L}_1(m)}([0,1]), \|\cdot \|_\infty)$ to ${\mathbb R}$ such
that $|\psi(x)|\leq C(1+\|x\|_\infty^2)$ for some positive constant $C$. Then 
$$
  \lim_{n \rightarrow \infty} {\mathbb E} \left ( \psi\left (\frac{S_n}{\sqrt n} \right) \right)= 
  {\mathbb E}((\psi(W)) \, .
$$
In particular
$$
\lim_{n \rightarrow \infty} {\mathbb E}  ( T_n )= 
  {\mathbb E}\left(\left (\max_{t \in [0,1]}\int |W_t(s)| \ m(ds)\right)^2 \right )  \, .
$$
\end{rem}

\noindent {\bf Proof of Theorem \ref{pifL1}.}
Note first that, in this adapted case,  all the conditions of Theorem \ref{thm:nonadapted}
are satisfied. Indeed, since for $m$ almost every $t$ the series 
(\ref{serieDR}) converge in ${\mathbb L}^1$, it follows that the series $\sum_{k=0}^\infty {\mathrm{Cov}}(Y_0(t), Y_k(t))$ converge, and then
$$
\lim_{n \rightarrow \infty} \frac{\|S_n(t)\|_2^2}{n} = {\mathrm{Var}}(Y_0(t)) + 2 \sum_{k=1}^\infty {\mathrm{Cov}}(Y_0(t), Y_k(t)) \, .
$$ 
Now, by definition of $L(t)$, 
$$
{\mathrm{Var}}(Y_0(t)) + 2 \sum_{k=1}^\infty {\mathrm{Cov}}(Y_0(t), Y_k(t)) \leq 2 L(t)\, .
$$
Hence the  condition \eref{condGordin} follows from  (\ref{condui}) and the fact that
$$
  C(t) \leq \sqrt{{\mathrm{Var}}(Y_0(t)) + 2 \sum_{k=1}^\infty {\mathrm{Cov}}(Y_0(t), Y_k(t))} \leq  \sqrt{2L(t)} \, .
$$
So,  the conclusion of Theorem \ref{thm:nonadapted} holds with the the covariance function defined by 
\eref{defcov}. 

As usual it remains to prove the tightness, which reduces through Ascoli's theorem to:
for any $\varepsilon>0$,
$$
\lim_{\delta \rightarrow 0} \limsup_{n \rightarrow \infty}
\frac 1 \delta
{\mathbb P}\left(\max_{1 \leq k \leq [n\delta]}  \int|S_k(t)| \ m(dt) > \sqrt n
\varepsilon \right)=0.
$$
But this follows straightforwardly from Lemma \ref{ui} below by applying Markov inequality at order 2.
The proof of Theorem \ref{pifL1} is complete. 
\begin{lem}\label{ui}
Assume that, for $m$-almost every $t$, the series 
defined in \eref{serieDR} converges in ${\mathbb L}^1$.
Assume moreover that the function $L$ defined in \eref{Lt} 
satisfies \eref{condui}.
Then the sequence $(T_n)_{n\geq 1}$ defined in \eref{maxui}
is uniformly integrable.
\end{lem}

\medskip

\noindent {\bf Proof of Lemma \ref{ui}.} We first note that, for any positive random variable $V$,
\begin{align*}
{\mathbb E} \left(\left ( \max_{1 \leq k \leq n} \int |S_k(t)| \ m(dt) \right)^2 V\right) &
\leq
{\mathbb E} \left(\left (  \int \sqrt V \max_{1 \leq k \leq n}  |S_k(t)|  \ m(dt) \right)^2\right)\nonumber\\
& \leq
\left (  \int \left \|\sqrt V \max_{1 \leq k \leq n}  |S_k(t)| \right \|_2 \ m(dt)
\right)^2\, .
\end{align*}
Taking $V={\bf 1}_{T_n>M}$, we obtain that 
\begin{equation}\label{majTn}
{\mathbb E}\left ( T_n {\bf 1}_{T_n>M}  \right) \leq
\frac 1 n \left (  \int \left \|  \left(\max_{1 \leq k \leq n}  |S_k(t)|\right) {\bf 1}_{T_n>M}  \right \|_2 \ m(dt)
\right)^2
\end{equation}
Applying Inequality (3.12) in \cite{DR} with $\lambda=0$, we get that
\begin{equation} \label{drineq1}
\left \|\max_{1 \leq k \leq n} |S_k(t)| \right \|_2^2 \leq 16 \sum_{k=1}^n
 \left \| Y_k(t) \sum_{i=k}^n {\mathbb E}_k(Y_i(t))
  \right \|_1 \leq 16 n L(t) \, .
\end{equation}
Using \eref{condui}, \eref{majTn},  \eref{drineq1} and the reverse Fatou Lemma, we infer that
$$
 \lim_{M \rightarrow \infty} \limsup_{n \rightarrow \infty}
 {\mathbb E}\left ( T_n {\bf 1}_{T_n>M}  \right)=0
$$
as soon as, for $m$-almost every $t$,
\begin{equation}\label{decadix}
\lim_{M \rightarrow \infty} \limsup_{n \rightarrow \infty}
\left \|  \frac {1}{\sqrt n}\left(\max_{1 \leq k \leq n}  |S_k(t)|\right) {\bf 1}_{T_n>M}  \right \|_2=0 \, . 
\end{equation}

It remains to prove \eref{decadix}. 
  In fact this follows quite easily from Proposition 1 in \cite{DR}. Indeed, since
  for $m$-almost every $t$, the series  defined in \eref{serieDR} converges in ${\mathbb L}^1$, it follows from this proposition that
 the sequence 
 $$
 \frac 1 n \left(\max_{1 \leq k \leq n}  |S_k(t)|\right)^2 
 $$
 is uniformly integrable for $m$-almost every $t$. 
 Hence \eref{decadix} holds as soon as 
 \begin{equation}\label{collaboratrice}
 \lim_{M \rightarrow \infty} \limsup_{n \rightarrow \infty}{\mathbb P}(T_n>M)=0 \, .
 \end{equation}
 Now, applying \eref{majTn} and \eref{drineq1}, 
 $$
 {\mathbb P}(T_n>M)\leq \frac{{\mathbb E}(T_n)}{M} \leq \frac {16}{M} \left( \int \sqrt{L(t)} \ m(dt) \right)^2 \, ,
 $$
 and \eref{collaboratrice} follows. This completes the proof of Lemma \ref{ui}.

\subsection{An invariance principle in ${\mathbb L}_1(m)$ for the empirical distribution function}\label{intro}
In this subsection, $S={\mathbb R}$, and $m$ is a $\sigma$-finite measure on ${\mathbb R}$ equipped
with  the Borel $\sigma$-field.
As in Section \ref{DefNot}, let $X_0$ be an ${\mathcal F}_0$-measurable and integrable real-valued random variable
with distribution function $F$.
Define the stationary sequence 
 ${\bf X}=(X_i)_{i \in {\mathbb Z}}$ by $X_i=X_0 \circ T^i$, and denote by $F_{X_k|{\mathcal F}_0}$ the conditional distribution function of $X_k$ given 
 ${\mathcal F}_0$. 
 
  The random variable $Y_k$
is then  defined by
$
Y_k(t)={\bf 1}_{X_k \leq t} -F(t) 
$,
in such a way that
$$
S_n= \sum_{k=1}^n Y_k = n(F_n-F) \, ,
$$
where $F_n$ is the empirical distribution function 
of $\{X_1, \ldots, X_n\}$. Note that $Y_0$ is a ${\mathbb L}^1(m)$-valued random variable as soon
as ${\mathbb E}(|X_0|) < \infty$. 

\begin{thm}\label{th:proj}
Assume that
\begin{equation}\label{proj}
\int \sqrt{ \sum_{k=0}^\infty \|F_{X_k|{\mathcal F}_0}(t)-F(t)\|_1 } \ m(dt) < \infty \, .
\end{equation}
Then $\{n^{-1/2} S_{[ns]}, s \in [0,1]\}$ converges in distribution
in  the space $D_{{\mathbb L}_1(m)}([0,1])$ to an ${\mathbb L}_1(m)$-valued Wiener process $W$. 
Moreover the explicit form of the covariance operator  
of $W$ is obtained {\it via} equation \eref{defcov} of Theorem \ref{thm:nonadapted}  by taking $Y_k(t)={\bf 1}_{X_k \leq t}-F(t)$.
\end{thm}
When applied to $\alpha$-dependent sequences as defined in Section \ref{DefNot},  Theorem 
\ref{th:proj} yields the following result. 

\begin{prop}\label{prop:alpha}
Let $B(t)=F(t)(1-F(t))$.
The condition
\begin{equation}\label{alpha}
\int \sqrt{ \sum_{k=0}^\infty \min\{\alpha_{1, {\bf X}}(k), B(t)\} } \ m(dt) < \infty
\end{equation}
implies the condition \eref{proj}, and hence the conclusion of Theorem
\ref{th:proj}. Moreover, the  covariance operator $\Lambda$ of $W$ can be expressed as follows: 
for any $f,g$ in ${\mathbb L}_\infty(m)$,
\begin{equation}\label{defcovbis}
\Lambda(f,g)= \sum_{k \in {\mathbb Z}} {\mathbb{E}}\left( \iint f(t)g(s) ({\mathbf{1}}_{X_0 \leq t}-F(t))   ({\mathbf{1}}_{X_k \leq s}-F(s)) \ m(dt) m(ds) \right)\, .
\end{equation}
\end{prop}

\noindent {\bf Proof of Theorem \ref{th:proj} and of Proposition \ref{prop:alpha}.}
Theorem \ref{th:proj} is a direct consequence of Corollary \ref{easycor} applied to the random variables
$
Y_k(t)={\bf 1}_{X_k \leq t} -F(t) 
$.
More precisely, since $|Y_0(t)|\leq 1$, the criterion \eref{easycrit} is exactly the criterion \eref{proj}.

It remains to prove Proposition \ref{prop:alpha}. We first quote that condition \eref{alpha} 
implies \eref{proj}: this follows easily from the two upper bounds given in \eref{b3}. 
It remains to prove that the covariance operator $\Lambda$ given in \eref{defcov} can be expressed 
as in \eref{defcovbis}. As usual, we identify a function $f$ in ${\mathbb L}^\infty(m)$ with an
element of the dual of ${\mathbb L}^1(m)$, and we write
$$
f(Y_k)= \int f(t) Y_k(t) \ m(dt).
$$
By Remark \ref{rem:extent}, we know that, for any $f$ in ${\mathbb L}^\infty(m)$,
\begin{equation}\label{firstlim}
\lim_{n \rightarrow \infty} \frac 1 n {\mathbb E}\left ( \left ( f(S_n) \right)^2 \right) = 
{\mathbb E}\left ( \left ( f(W_1) \right)^2 \right) = \Lambda(f,f) \, .
\end{equation}
Now, if we can prove that, for any $f, g$ in ${\mathbb L}^\infty(m)$,
\begin{equation}\label{cov}
\sum_{k\in {\mathbb Z}} |{\mathrm{Cov}}(f(Y_0), g(Y_k))| < \infty \, ,
\end{equation}
then the series
\begin{equation*}
  \bar \Lambda(f,g)= \sum_{k\in {\mathbb Z}} {\mathrm{Cov}}(f(Y_0), g(Y_k)) 
\end{equation*}
is well defined, and 
\begin{equation}\label{secondlim}
\lim_{n \rightarrow \infty} \frac 1 n {\mathbb E}\left ( \left ( f(S_n) \right)^2 \right) = 
 \bar \Lambda(f,f) \, .
\end{equation}
From \eref{firstlim} and \eref{secondlim}, we infer that, for any $f$ in ${\mathbb L}^\infty(m)$,
$
  \Lambda(f,f)= \bar \Lambda (f,f) 
$.
Applying this equality to $f$, $g$, and $f+g$ it follows that, for any $f, g$ in ${\mathbb L}^\infty(m)$,
$$
\Lambda(f,g)= \bar \Lambda (f,g) \, ,
$$
which is the desired result. To prove \eref{cov}, we first note that
$$
|{\mathrm{Cov}}(f(Y_0), g(Y_k))| \leq \|f\|_\infty \| g\|_\infty
\int \int \|({\mathbf{1}}_{X_0 \leq t}-F(t))  {\mathbb E}_0 ({\mathbf{1}}_{X_k \leq s}-F(s)) \|_1 \ m(dt)m(ds) \, .
$$
Now 
$$
\|({\mathbf{1}}_{X_0 \leq t}-F(t))  {\mathbb E}_0 ({\mathbf{1}}_{X_k \leq s}-F(s)) \|_1
\leq \min \{\alpha_{1, \bf X}(k) ,  2B(t) , 2 B(s)\} \, .
$$
Hence
$$
\sum_{k=0} \|({\mathbf{1}}_{X_0 \leq t}-F(t))  {\mathbb E}_0 ({\mathbf{1}}_{X_k \leq s}-F(s)) \|_1
\leq 2 \sqrt{ \sum_{k=0}^\infty \min\{\alpha_{1, \bf X}(k), B(t)\} } \sqrt{ \sum_{k=0}^\infty \min\{\alpha_{1, \bf X}(k), B(s)\} }
\, .
$$
This implies  that
$$
\sum_{k\in {\mathbb Z}} |{\mathrm{Cov}}(f(Y_0), g(Y_k))|
\leq
2 \|f\|_\infty \|g\|_\infty \left ( \int \sqrt{ \sum_{k=0}^\infty \min\{\alpha_{1, \bf X}(k), B(t)\} } \ m(dt)
\right )^2 \, ,
$$
and \eref{cov} follows from   \eref{alpha}. This completes the proof of 
Proposition \ref{prop:alpha}.

\section{Quantile conditions} \label{quantile}
\setcounter{equation}{0}
As a consequence of the results by D\'ed\'e \cite{D} or  Cuny \cite{C} (see the condition 
\eref{conseq} of Subsection \ref{sec:CD}) we know that
the conclusion of Theorem \ref{th:proj} holds as soon as
\begin{equation}\label{CD}
\sum_{k=0}^\infty \frac{1}{\sqrt {k +1}} \int  \|F_{X_k|{\mathcal F}_0}(t)-F(t)\|_2  \ m(dt) < \infty \, .
\end{equation}
Moreover, it follows from \cite{C} that the condition \eref{CD} also implies the strong invariance principle.

Let $B(t)=F(t)(1-F(t))$. 
As quoted by D\'ed\'e (2009), the condition \eref{CD} is implied by
\begin{equation}\label{SophieCond}
\sum_{k=0}^\infty \frac{1}{\sqrt{k +1}} \int \sqrt {\min \{ \alpha_{1, {\bf X}}(k), B(t) \}}  \ m(dt) \, .
\end{equation}
The conditions \eref{alpha} of Proposition \ref{prop:alpha} and the condition \eref{SophieCond} are not easy to compare. However, if either
$m$ has finite mass or $X_0$ is bounded, then \eref{alpha}
is equivalent to
\begin{equation}
\sum_{k=1}^\infty \alpha_{1, {\bf X}}(k) < \infty
\end{equation}
and \eref{SophieCond} is equivalent to
\begin{equation}
\sum_{k=1}^\infty \sqrt{\frac{\alpha_{1, {\bf X}}(k)}{k}} < \infty\, .
\end{equation}
Hence, in that case, the condition \eref{alpha} is weaker than the condition \eref{SophieCond},
and is in fact equivalent to the minimal condition to get the central limit theorem for partial sums
of stationary $\alpha$-dependent sequences of bounded random variables.

We shall now focus on the
the case where $m=\lambda$ is the Lebesgue measure on ${\mathbb R}$. In that case, the condition \eref{SophieCond} is equivalent to
\begin{equation}\label{SophieCondbis}
\sum_{k=0}^\infty \frac{1}{\sqrt {k+1}} \int_0^\infty \sqrt {\min \{ \alpha_{1, {\bf X}}(k), H(t) \}}  \ dt \, .
\end{equation}
and the condition \eref{alpha} is equivalent to \eref{condalpha}.
We shall see that
the condition \eref{condalpha} is always weaker than the condition \eref{SophieCondbis}.
The first step is to express \eref{condalpha} and  \eref{SophieCondbis} in terms of the quantile
function of $X_0$, as done in \cite{DMR} for the invariance principle of stationary $\alpha$-mixing sequences. More precisely, we shall
compare the three following conditions:
\begin{align}
  \int_0^1 \alpha^{-1}(u) Q^2(u) du & < \infty\, , \label{DMR} \\
  \int_0^1 \frac{ \alpha^{-1}(u) Q(u)} {\sqrt {\int_0^u \alpha^{-1}(x) \ dx }} du &< \infty \, , \label{DM}\\
  \int_0^1 \frac{ \sqrt{ \alpha^{-1}(u)} Q(u)}{\sqrt u} du & < \infty \label{D}\, , 
\end{align}
where for simplicity we denote by $\alpha^{-1}$ the function
$\alpha_{1,{\bf X}}^{-1}$  defined in \eref{alphamoins1}. 
The condition \eref{DMR} has been introduced by Doukhan, Massart and Rio \cite{DMR}, but in that paper
the function $\alpha^{-1}$ is defined with the
$\alpha$-mixing coefficients  of Rosenblatt \cite{R}. These authors showed that \eref{DMR}
 implies the functional central limit theorem for the Donsker line
$$
\left \{ \frac{1}{\sqrt n} \sum_{k=1}^{[nt]} \big (X_k- {\mathbb E}(X_k)), \  t \in [0,1] \right \} \, ,
$$
and that it is optimal in a precise sense.
The optimality of this condition has been further discussed in a
paper by Bradley \cite{B}. The fact that, for ergodic sequences, this functional central limit theorem remains true 
with the much weaker coefficients $\alpha_{1,{\bf X}}(k)$ is a consequence of a result by 
Dedecker and Rio \cite{DR}. 

Concerning these three quantile conditions, our first result is Proposition \ref{equiv} below.
\begin{prop} \label{equiv} The following equivalences hold
\begin{enumerate}
\item The condition \eref{DMR} is equivalent to
\begin{equation}\label{DMRbis}
\int_0^{\infty} t \left( \sum_{k=0}^\infty \min\{\alpha_{1,{\bf X}}(k), H(t)\} \right ) \ dt < \infty \, .
\end{equation}
\item The condition \eref{DM} is equivalent to \eref{condalpha}.
\item
The condition \eref{D} is equivalent to \eref{SophieCondbis}.
\end{enumerate}
\end{prop}
The hierarchy of these quantile conditions is given in Proposition \ref{hier} below.
\begin{prop}\label{hier}
The following implications hold:
\eref{D} $\Rightarrow$ \eref{DM} $\Rightarrow$ \eref{DMR}.
\end{prop}

\begin{rem}
At this point, it should be noticed that these three conditions are in fact very close.
Indeed, by a simple application of Cauchy-Schwarz inequality, for any $b>1/2$,
$$
  \eref{DMR} \Rightarrow \int_0^1 \frac{ \sqrt{ \alpha^{-1}(u)} Q(u)}{\sqrt u |1+\ln(u)|^b} du  < \infty
$$
and the condition on right hand is a  slight reinforcement of \eref{D}.
\end{rem}

\noindent{\bf Proof of Proposition \ref{equiv}.} Assume that $\sum_{k\geq 0} \alpha_{1,{\bf X}}(k)
 < \infty$. Then the function
$S$ defined
on ${\mathbb R}^+$ by
\begin{equation}\label{def:S}
S(t)= \sum_{k=0}^\infty \min\{\alpha_{1,{\bf X}}(k), H(t)\}
\end{equation}
is finite and  non-increasing.

\medskip

\noindent{\it Proof of Item 1.}  By a simple change of variables,  we see that the condition \eref{DMRbis} is equivalent
to
$$
\int_0^\infty S \left( \sqrt t\right ) \ dt < \infty \, .
$$
Since
\begin{equation}\label{S}
S(t)= \sum_{k=0}^\infty \int_0^1 {\bf 1}_{u \leq \min\{\alpha_{1,{\bf X}}(k), H(t)\}} \ du=
\int_0^{H(t)} \alpha^{-1}(u) du \, ,
\end{equation}
it follows that
\begin{align*}
\int_0^\infty S \left(\sqrt t \right) \ dt & = \int_0^\infty  \left(\int_0^1 \alpha^{-1}(u) {\bf 1}_{u \leq H\left(\sqrt t\right)}  \ du \right) \ dt \\
                               & = \int_0^1 \alpha^{-1}(u) 
                               \left(\int_0^\infty {\bf 1}_{t \leq Q^2(u)}  \ dt \right) \ du \\
                               & = \int_0^1  \alpha^{-1}(u) Q^2(u)  \ du \, ,
\end{align*}
which concludes the proof of Item 1.

\medskip

\noindent {\it Proof of Item 2.}  Starting from \eref{S},
it follows that
\begin{equation}\label{un}
\int_0^{\infty} \sqrt{ \sum_{k=0}^\infty \min\{\alpha_{1,{\bf X}}(k), H(t)\} } \ dt
= \int_0^{\infty} \sqrt{ \int_0^{H(t)} \alpha^{-1}(u) \ du} \ dt \, .
\end{equation}
Let
$$
G_\alpha(x) =\sqrt{\int_0^x \alpha^{-1}(u) \ du } \, .
$$
From \eref{un}, we infer that
\begin{align*}
\int_0^{\infty} \sqrt{ \sum_{k=0}^\infty \min\{\alpha_{1,{\bf X}}(k), H(t)\} } \ dt
&= \int_0^\infty \int_0^1 {\bf 1}_{v \leq G_\alpha(H(t))} \ dv \ dt \\
&= \int_0^1 \int_0^\infty {\bf 1}_{t \leq Q ( G^{-1}_\alpha(v))} \ dt \ dv \\
&= \int_0^1 Q \circ G_{\alpha}^{-1} (v)  \ dv \, .
\end{align*}
Making the change of variables $u=G_\alpha^{-1}(v)$, the result follows.

\medskip

\noindent {\it Proof of Item 3.}
Note first that
\begin{align}
 \int_0^\infty \sqrt {\min \{ \alpha_{1,{\bf X}}(k), H(t) \}}  \ dt & = \int_0^\infty 
 \left( \int_0^1 {\bf 1}_{u^2 \leq \alpha_{1,{\bf X}}(k)} {\bf 1}_{u^2 \leq H(t)} \ du \right) \ dt \\
                                                        & =  \int_0^1 Q(u^2) {\bf 1}_{u^2 \leq 
                                                        \alpha_{1,{\bf X}}(k)} \ du \, .
\end{align}
Now $u^2 \leq \alpha_{1,{\bf X}}(k)$ if and only if $k \leq \alpha^{-1}(u^2)$. Hence, there exists two positive constants $A$ and $B$  such that
$$
A \sqrt{\alpha^{-1}(u^2)} \leq  \sum_{k=1}^\infty \frac{1}{\sqrt k} {\bf 1}_{u^2 \leq \alpha_{1,{\bf X}}(k)} \leq B \sqrt{\alpha^{-1}(u^2)} \, .
$$
Finally
$$
\sum_{k=1}^\infty \frac{1}{\sqrt k} \int_0^\infty \sqrt {\min \{ \alpha_{1,{\bf X}}(k), H(t) \}}   \ dt
<  \infty \quad \text{iff} \quad  \leq  \int_0^1 \sqrt{\alpha^{-1}(u^2)} Q(u^2) du < \infty \, .
$$
Making the change of variables $v=u^2$, the result follows.

\medskip

\noindent {\bf Proof of Proposition \ref{hier}.}
 Since the function $\alpha^{-1}$ is non-increasing,
one has
$$
 \int_0^u \alpha^{-1}(x) \ dx \geq u \alpha^{-1}(u) \, ,
$$
which proves that \eref{D} implies \eref{DM}.

It remains to prove that \eref{DM} implies \eref{DMR}. By Proposition \ref{equiv}, it is equivalent to prove that \eref{condalpha} implies \eref{DMRbis}. If \eref{condalpha} holds, then the function $S$ defined
on ${\mathbb R}^+$ by
\eref{def:S}
is finite and  non-increasing. Hence, using again \eref{condalpha},
$$
t \sqrt{S(t)} \leq 2 \int_{t/2}^t \sqrt{S(s)} \ ds \leq C\, , \quad \text{with} \quad
C= 2 \int_{0}^\infty \sqrt{S(s)} \ ds \, .
$$
Consequently
$
t S(t) \leq C \sqrt{S(t)}
$,
proving that \eref{condalpha} implies \eref{DMRbis}.

\subsection{Sufficient conditions}
In this subsection, we give some simple conditions on $\alpha_{1,{\bf X}}(k)$ and $H$ under which 
(\ref{DM}) (and hence \eref{condalpha}) is satisfied. 
\begin{prop}\label{prop:suff}
The following conditions imply \eref{DM}:
\begin{enumerate}
\item
$${\mathbb E}(|X_0|^p) < \infty  \ \text{for some $p>2$, and} \quad
   \sum_{k>0} \frac{\left(\alpha_{1,{\bf X}}(k)\right)^{\frac{p-2}{2(p-1)}}} {k^{\frac{p-2}{2(p-1)}}}< \infty \, .
$$
\item
$$
H(t) =O(t^{-p}) \  \text{for  some $p>2$, and} \quad
 \sum_{k>0} \frac{\left(\alpha_{1,{\bf X}}(k)\right)^{\frac{p-2}{2p}}} {\sqrt k}< \infty \, .
$$
\item
$$
  \int_0^\infty \big( H(t)\big)^{\frac{a-1}{2a}} \ dt < \infty
  \quad \text{and $\alpha_{1,{\bf X}}(k) =O\Big( \frac{1}{k^a}\Big)$ for some $a>1$.}
$$
\item
$$
  \int_0^\infty \left ( \ln \left ( 1 +  \frac{1}{H(t)}\right)\right)^{-\frac{(a-1)}{2}} dt < \infty
  \quad \text{and $\alpha_{1,{\bf X}}(k) =O\Big( \frac{1}{k (\ln(k))^a}\Big)$ for some $a>1$.}
$$
\item
$$
  \int_0^\infty \sqrt{ H(t) \big | \ln (H(t))\big |} \ dt < \infty
  \quad \text{and $\alpha_{1,{\bf X}}(k) =O(a^k)$ for some $a<1$.}
$$
\end{enumerate}
\end{prop}

\noindent{\bf Proof of Proposition \ref{prop:suff}.}

\noindent{\it Proof of Item 1.}
Since \eref{D} implies \eref{DM}, it suffices to prove  that Item 1 implies \eref{D}.
Applying Cauchy Schwarz, we obtain that
$$
\int_0^1 \frac{ \sqrt{ \alpha^{-1}(u)} Q(u)}{\sqrt u} du \leq \left (\int_0^1 Q(u)^p du\right)^{\frac 1 p}
\left( \int_0^1 \left(\frac{\alpha^{-1}(u)}{u}\right )^{\frac{p}{2(p-1)}}
\ du \right )^{\frac{p-1}{p}} \, .
$$
Since ${\mathbb E}(|X_0|^p) < \infty$, the first integral on right hand is finite. It remains to prove
that
$$
\int_0^1 \left(\frac{\alpha^{-1}(u)}{u}\right )^{\frac{p}{2(p-1)}}
\ du < \infty \, .
$$
By definition of $\alpha^{-1}$, this is equivalent to
$$
  \sum_{k>0} k^{\frac{p}{2(p-1)}} \int_{\alpha_{1,{\bf X}}(k+1)}^{\alpha_{1,{\bf X}}(k)}
  u^{-\frac{p}{2(p-1)}} du  < \infty \, .
$$
The last condition means exactly that
$$
\sum_{k>0} k^{\frac{p}{2(p-1)}}\left( \left(\alpha_{1,{\bf X}}(k)\right)^{\frac{p-2}{2(p-1)}}- 
\left(\alpha_{1,{\bf X}}(k+1)\right)^{\frac{p-2}{2(p-1)}}
\right) < \infty \, ,
$$
which is equivalent to the condition of Item 1.

{\it Proof of Item 2.} Again, it suffices to prove that Item 2 implies \eref{D}.
Now, the condition $H(t)=O(t^{-p})$ is equivalent to $Q(u)=O(u^{-1/p})$. Hence, the condition
\eref{D} holds as soon as
$$
  \int_0^1 \frac{\sqrt{\alpha^{-1}(u)}}{u^{\frac 1 p + \frac 1 2}} du < \infty \, .
$$
By definition of $\alpha^{-1}$, the last condition means exactly that
$$
\sum_{k>0} \sqrt{k} \left( \left(\alpha_{1,{\bf X}}(k)\right)^{\frac{p-2}{2p}}- 
\left(\alpha_{1,{\bf X}}(k+1)\right)^{\frac{p-2}{2p}}
\right) < \infty \, ,
$$
which is equivalent to the condition of Item 2.

{\it Proofs of Item 3, 4 and 5.}
For the proof of these points, we start from condition \eref{condalpha} which is equivalent to
\eref{DM}. Since we can control the behavior of $\alpha_{1,{\bf X}}(k)$, we can give upper bounds
for the function $S$ defined by \eref{def:S}.

If $\alpha_{1,{\bf X}}(k) =O\left( \frac{1}{k^a}\right)$ for some $a>1$, then $S(t)=O\left (\left( H(t)\right)^{\frac{a-1}{a}}\right )$.

If $\alpha_{1,{\bf X}}(k) =O\left( \frac{1}{k (\ln(k))^a}\right)$ for some $a>1$, then
$$
S(t)=O \left (
\left ( \ln \left ( 1 +  \frac{1}{H(t)}\right)\right)^{-\frac{(a-1)}{2}}
\right)
\, .
$$

If $\alpha_{1,{\bf X}}(k) =O(a^k)$ for some $a<1$, then $S(t)=O\left(H(t) \left | \ln (H(t))\right 
|\right)$.

Item 3, 4 and 5 follow from these upper bounds and condition \eref{condalpha}.

\end{document}